\documentclass[12pt]{amsart}
\usepackage{amsfonts}
\usepackage{amssymb}
\usepackage[dvipsnames]{xcolor}
\usepackage[pdftex]{graphicx}
\date{July 2019}

\renewcommand{\epsilon}{\varepsilon}
\renewcommand{\phi}{\varphi}

\newtheorem{Lemma}{Lemma}[section]
\newtheorem{Theorem}{Theorem}[section]

\newtheorem{Corollary}{Corollary}[section]
\newtheorem{Definition}{Definition}[section]

\newtheorem{Remark}{Remark}[section]

\newtheorem{Question}{Question}[section]

\begin{document}

\address{Department of Mathematical Sciences, Binghamton University, 4400 Vestal Parkway East
Binghamton, New York 13902-6000, USA.}
\author{Alexander Borisov}
\thanks{}
\title[]{Frameworks for two-dimensional Keller maps}

\begin{abstract} A Keller map is a counterexample to the Jacobian Conjecture. In dimension two every such map, if exists, leads to a complicated set of conditions on the map between the Picard groups of suitable compactifications of the affine plane. This is essentially a combinatorial problem. Several solutions to it (``frameworks") are described in detail. Each framework corresponds to a large system of equations, whose solution would lead to a Keller map. 
\end{abstract}

\email{borisov@math.binghamton.edu}

\maketitle

\section{Introduction} 

Suppose $f(x,y)$ and  $g(x,y)$ are two polynomials with complex coefficients. The classical Jacobian Conjecture (due to Keller, \cite{Keller}) asserts the following.

{\bf Conjecture.} (Jacobian Conjecture in dimension two) If the Jacobian of the pair of two-variable complex polynomials $(f,g)$ is a non-zero constant, then the map $(x,y)\mapsto (f(x,y),g(x,y))$ is invertible. Note that the opposite is clearly true, because the Jacobian of any polynomial map is a polynomial, and, when the map is invertible, it must have no zeros, so it is a constant.

This conjecture has a long history. See an excellent survey of Miyanishi \cite{MiyanishiAAG} for some references.  For the more algebraic approaches see the survey of van den Essen \cite{Essen}. The term ``Jacobian Conjecture"  was coined by Abhyankar (cf.  \cite{AbhyankarTataLectures}).

The approach of this paper is based on the birational geometry of complete surfaces and combinatorial properties of the graphs of curves at infinity, as in the papers \cite{AmpleRamification} and \cite{DivisorialValuations}. All varieties are over $\mathbb C.$

From the point of view of a birational geometer, the most natural approach to the two-dimensional Jacobian Conjecture is the following. Suppose a counterexample, called a Keller map, exists. It gives a rational map from $X={\mathbb P}^2$ to $Y={\mathbb P}^2.$ After a sequence of blowups of  outside of $\mathbb A^2$, we can get a surface $Z$ with two morphisms: $\pi : Z \to X$ (projection onto the source $\mathbb P^2$) and $\phi : Z \to Y$ (the lift of the original rational map).

Note that $Z$ contains a Zariski open subset isomorphic to $\mathbb A^2$, and that its complement, $\pi ^* ((\infty))$, is a tree of smooth rational curves. We will call these curves exceptional, or curves at infinity. The structure of this tree is easy to understand inductively, as it is built from a single curve $(\infty)$ on ${\mathbb P}^2$ by a sequence of two operations: blowing up a point on one of the curves or blowing up a point of intersection of two curves. It is important to note that the exceptional curves on $Z$ may behave very differently with respect to the map $\phi$. More precisely, there are four types of curves $E$.

type 1) $\phi (E)$ is a curve, $\phi (E) \cap \mathbb A^2 = \emptyset$ (i.e., $\phi (E)=(\infty)$)
 
type 2) $\phi (E)$ is a point, not in $\mathbb A^2$

type 3) $\phi (E)$ is a curve, $\phi (E) \cap \mathbb A^2 \neq \emptyset$ (i.e., $\phi (E) \neq (\infty)$)

type 4) $\phi (E)$ is a point in $\mathbb A^2$

The curves of type $3$ are known as di-critical components (cf., e.g. \cite{Abhyankar2010}, \cite{Orevkov}). A trivial topological argument implies that such curve must exist in any counterexample to the Jacobian Conjecture. A slightly stronger result, that the map must be ramified in at least one di-critical component, was proved (possibly, reproved) in \cite{AmpleRamification}, Theorem 3.1.

As was already done before, in particular by Orevkov, we are going to apply a sequence of blowups at infinity to the target surface, in the attempt to ``get a closer view" of the Keller map. After adjusting the target, we again resolve the map. Slightly abusing the notation, we will call the new target surface \(Y\) and the new resolution surface \(Z\).  In this new situation, the classification of the exceptional curves on the new surface \(Z\)   into four types still makes sense, just skip the parts in parentheses for types 1 and 3. Note that some of the curves that were classified as type 2 when the target was $\mathbb P^2$ may now be of type \(1\). No other type changes can occur. 

We can consider the Stein factorization of  the morphism $\phi : Z \to Y.$ That is, we factor it into a composition of two morphisms, birational and finite: $Z \longrightarrow W \longrightarrow Y$.  Here the first morphism is birational and will be denoted  by $\tau$, and the second one is finite and will be denoted by $\rho.$ The surface $W$ is  algebraic and normal. In what follows, we will use the intersection theory for complete normal surfaces due to Mumford. Suppose $K_{W}$ is the canonical class of $W$, as the Weil divisor class modulo numerical equivalence. The augmented canonical class is, by definition, $\bar{K}_{W}=K_{W}+\sum E_i$, where \(E_i\) are the images of all exceptional curves of types 1 and 3 (cf. \cite{Kollar_pairs}). All the curves of types 2 and 4 are contracted by $\tau$.

{\bf Keller Map Adjunction Formula}
The two conditions of being a Keller map, $\mathbb A^2$ is mapped to $\mathbb A^2$ and no ramification on $\mathbb A^2,$ can be combined in one formula in the Mumford Picard group on $W$ (cf. \cite{AmpleRamification}):
$$\bar{K}_{W}=\rho^*(\bar{K}_Y)+\bar R,$$
where $\bar{R}=\sum \limits_{type(E_i) = 3} e_i E_i,$ where $e_i$ is the ramification index of $\phi$ at $E_i$.

It was proven in \cite{AmpleRamification} that when $Y=\mathbb P^2$ the curve $\pi^{-1}_*(\infty)$ on $Z$ is of type 2. Moreover, the surface $W$ has a rather simple structure, apart from one  possibly complicated point, $\tau(\pi^{-1}_*(\infty))$. Therefore, it makes sense to start modifying the target, by blowing up $\phi(\pi^{-1}_*(\infty))$, until it becomes a curve.  There are several restrictions on this process, some more complicated than the others. In particular, on $Y$ all $\bar{K}$ labels are non-positive, and all determinant labels are positive (cf. \cite{DivisorialValuations} for the definitions).  This is how all our frameworls were constructed, by hand. It should be stressed that the existence of these frameworks is no miracle: the obstructions seem to be of the ``inequality" type rather than ``congruence" type, or anything trickier.  One should expect  infinitely many frameworks similar to the ones presented in this paper.

The paper is organized as follows. In section 2 some preliminaries are discussed and the first framework is constructed. Section 3 discusses it in further details, as well as the attempts to construct a Keller map based on it. Section 4 contains the second framework. Section 5 describes some frameworks that are related to the first framework, and a more complicated framework. Section 6 discusses the origins of these frameworks and some natural next steps for  attacking the two-dimensional Jacobian Conjecture from this direction.  

\section{Preliminaries, Notation, and First Framework}

Throughout the paper, we will be dealing primarily with smooth compactifications of $\mathbb A^2$, obtained from $\mathbb P^2$ by a sequence of blowups of  points  and contractions of $(-1)$-curves outside of $\mathbb A^2$. To every such surface we can associate a graph of curves ``at infinity" (i.e. outside of $\mathbb A^2$). The vertices of this graph are the ``curves at infinity", i.e. the irreducible components of the complement of $\mathbb A^2$. The two vertices are joined by an edge whenever the two curves intersect. Here we assume that the curves are in simple normal crossing. This will be automatically achieved, as long as we never contract a $(-1)$-curve that intersects three or more other curves at infinity. Note also that this graph is a tree.

Because every divisor on $\mathbb A^2$ is a divisor of a function, the classes of curves at infinity generate the Picard group of our compactified surface. Moreover, they form its basis (the fact that can be easily proven by induction). Therefore the graph of exceptional curves together with the self-intersection numbers of the curves completely determine the structure of the intersection form on our surface. Note that the self-intersection numbers change under blowups and contractions, so they are not the invariants of the divisorial valuations defined by the curves at infinity. To get around that, two other labels for these curves were introduced in \cite{AmpleRamification} and \cite{DivisorialValuations}: the $\bar{K}$ label and the determinant label. These labels are invariant under polynomial automorphisms of $\mathbb A^2$. Modulo that, the valuations with given labels  form finitely many  families (cf.  \cite{DivisorialValuations}). We will not use the determinant labels until section 5, and will define and discuss them there. Here is the definition of the $\bar{K}$ labels, that will be used a lot throughout the paper. 

\begin{Definition} The $\bar{K}$ label of a curve at infinity is the coefficient in the expansion of $\bar{K}=K+\sum E_i$ in the basis $\{E_i\}$ of the Picard group of our surface. Here $K$ is the canonical class, and the $\bar{K}$ is the augmented canonical class, which is the sum of $K$ and the ``boundary", that is, naturally, the sum of all curves at infinity.
\end{Definition}

The $\bar{K}$ labels are easier to work with than the self-intersection labels, because once a curve is created its $\bar{K}$ label no longer changes. When a curve is created by blowing up a point on one of the curves at infinity, its $\bar{K}$ label is obtained by adding $1$ to the $\bar{K}$ label of its ``parent" curve. If it is obtained by blowing up the point of intersection of two curves, its $\bar{K}$ label is simply the sum of the $\bar{K}$ labels of its two ``parents". On the original $\mathbb P^2$ the line at infinity has $\bar{K}$ label $(-2).$

The following observation is easy but very important.

\begin{Lemma} Suppose $\phi: Z\to Y$ is the Keller map (in the notation of the Introduction), and $E_i$ is a type 1 curve at infinity on $Z.$ Suppose that $\phi(E_i)$ is $F_i$ and the ramification index of $\phi $ at $E_i$ is $e_i$. Then the $\bar{K}$ label of $E_i$ is $e_i$ times the $\bar{K}$ label of  $F_i$.
\end{Lemma}

\proof
This follows directly from the formula $\bar{K}_{W}=\rho^*(\bar{K}_Y)+\bar R$.
\endproof

As was noted in \cite{AmpleRamification}, for every curve at infinity $E$ with non-zero $\bar{K}$ label one can recover its self-intersection index from its $\bar{K}$ label and  the $\bar{K}$ labels of its neighbors. Indeed, suppose $E$ has $k$ neighbors, $E_1,...,E_k$ with $\bar{K}$ labels  $a_i$, and the $\bar{K}$ label of $E$ is $a$. By the adjunction formula for $E,$
$$-2=(K+E)\cdot E=\bar{K}\cdot E -k=aE^2+\sum \limits_{i=1}^{k}a_i -k,$$
$$-E^2=\frac{1}{a}(\sum \limits_{i=1}^{k}a_i-k+2)$$
In particular, when $k=2$, we get $-E^2=\frac{a_1+a_2}{a}$.

On the other hand, for curves with $\bar{K}$ label $0$ it is not possible to recover $E^2$ from the $\bar{K}$ labels, and this is significant. In general, the curves with $\bar{K}$ label $0$ seem to play an important and somewhat mysterious role in the subject (see section 6 for more on that). As a result, in our figures we will note the $\bar{K}$ labels for all curves (under the corresponding vertices of the graph) and the self-intersections of the curves with $\bar{K}$ label $0$ (in parentheses) under their zero $\bar{K}$ labels.

Now we are ready to construct the first framework. We will start with the construction of the target surface $Y$, because it is easier, then construct the source surface $Z$, and then write down $\phi_*$ and $\phi^*$.

Starting with the projective plane, we first blow up a point on the line at infinity, then the intersection of the newly created exceptional curve with the strict pullback of the line at infinity, then again the intersection of the newly created curve and the pullback of the line at infinity. As a result, we get a surface with the following graph:

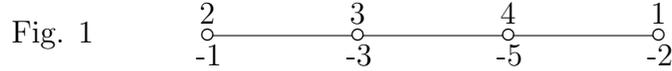
\setlength{\unitlength}{1cm}
\begin{picture}(10,1.5)
\put(0.5,0.5){Fig. 1}
\multiput(3,0.5)(2,0){4}{$\circ$}
\multiput(3.17,0.595)(2,0){3}{\line(1,0){1.86}}
\put(2.95,0.2){-$1$}
\put(4.95,0.2){-$3$}
\put(6.95,0.2){-$5$}
\put(8.95,0.2){-$2$}
\put(3,0.75){$2$}
\put(5,0.75){$3$}
\put(7,0.75){$4$}
\put(9,0.75){$1$}
\end{picture}

Here the integers below the vertices are the $\bar{K}$ labels, and the natural numbers above them simply stand for the order in which the curves were constructed.

Now we blow up a point on the curve that was last constructed, then a point on the newly constructed curve, and again, and again, until we get to the curve with the $\bar{K}$ label $0$:

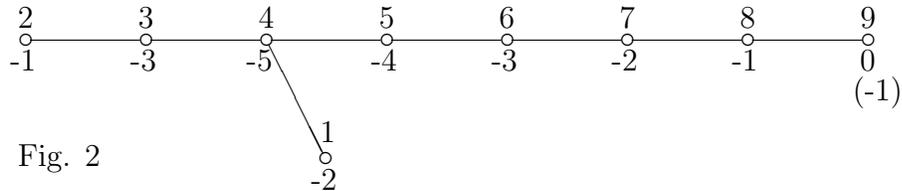
\setlength{\unitlength}{.8cm}
\begin{picture}(10,1.5)
\put(0.5,-1.5){Fig. 2}
\multiput(0.5,0.48)(2,0){8}{$\circ$} 
\multiput(0.7,0.6)(2,0){7}{\line(1,0){1.83}}
\put(0.37,0.1){-$1$}
\put(2.37,0.1){-$3$}
\put(4.3,0.1){-$5$}
\put(6.37,0.1){-$4$}
\put(8.37,0.1){-$3$}
\put(10.37,0.1){-$2$}
\put(12.37,0.1){-$1$}
\put(14.5,0.1){$0$}
\put(14.37,-0.4){(-$1$)}

\put(0.5,0.8){$2$}
\put(2.5,0.8){$3$}
\put(4.5,0.8){$4$}
\put(6.5,0.8){$5$}
\put(8.5,0.8){$6$}
\put(10.5,0.8){$7$}
\put(12.5,0.8){$8$}
\put(14.5,0.8){$9$}

\put(4.67,0.52){\line(1,-2){0.9}}
\put(5.48,-1.49){$\circ$} 
\put(5.37,-1.89){-$2$}
\put(5.52,-1.1){$1$}
\end{picture}
\vskip 1.8cm
As mentioned before, the number inside the parentheses indicates the self-intersection of the curve with the $\bar{K}$ label $0$.

Next, we blow up the point of intersection of the last two curves, and the point of intersection of the two curves with $\bar{K}$ labels (-1). Finally, we blow up a point on the newly created curve to get another curve with $\bar{K}$ label (-1), and a point on that curve, to get our surface $Y$:

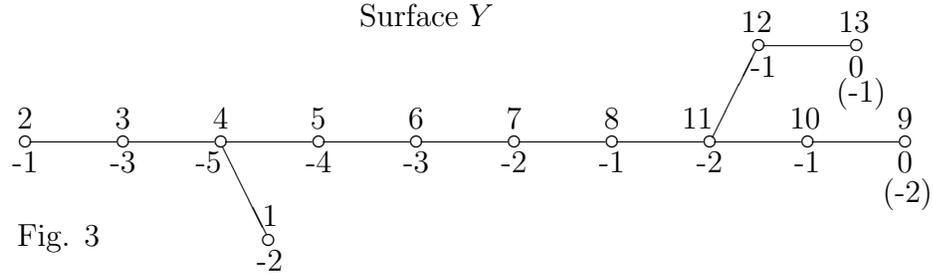
\setlength{\unitlength}{.65cm}
\begin{picture}(10,3.5)
\put(7.5,3){Surface $Y$}
\put(0.5,-1.5){Fig. 3}
\multiput(0.5,0.48)(2,0){10}{$\circ$} 
\multiput(0.75,0.63)(2,0){9}{\line(1,0){1.78}}
\put(0.4,0){-$1$}
\put(2.4,0){-$3$}
\put(4.15,0){-$5$}
\put(6.4,0){-$4$}
\put(8.4,0){-$3$}
\put(10.4,0){-$2$}
\put(12.4,0){-$1$}
\put(14.4,0){-$2$}
\put(16.4,0){-$1$}
\put(18.5,0){$0$}
\put(18.2,-0.6){(-$2$)}

\put(0.5,0.9){$2$}
\put(2.5,0.9){$3$}
\put(4.5,0.9){$4$}
\put(6.5,0.9){$5$}
\put(8.5,0.9){$6$}
\put(10.5,0.9){$7$}
\put(12.5,0.9){$8$}
\put(14.1,0.9){$11$}
\put(16.3,0.9){$10$}
\put(18.5,0.9){$9$}

\put(4.7,0.52){\line(1,-2){0.89}}
\put(5.47,-1.53){$\circ$} 
\put(5.4,-2){-$2$}
\put(5.51,-1.1){$1$}

\put(14.72,0.72){\line(1,2){0.9}}
\put(15.5,2.45){$\circ$} 
\put(15.75,2.6){\line(1,0){1.78}}
\put(17.5,2.45){$\circ$}
\put(15.3,2.87){$12$}
\put(17.3,2.87){$13$}
\put(15.5,1.97){-$1$}
\put(17.5,1.97){$0$}
\put(17.25,1.45){(-$1$)}
\end{picture}
\vskip 1.8cm

Now we will construct the surface $Z,$ which is considerably more complicated. We start with $\mathbb P^2$ and blow up a point on the line at infinity, and then a point on the last curve, to get a curve with $\bar{K}$ label $0$. Then we blow up the intersection of the last two curves and then the intersection of the two curves with $\bar{K}$ labels (-1), to get a curve with $\bar{K}$ label $(-2).$ Here is the resulting graph:

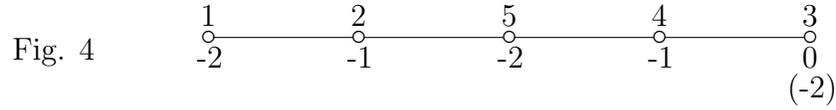
\setlength{\unitlength}{1cm}
\begin{picture}(10,1.5)
\put(0.5,0.3){Fig.  4}
\multiput(3,0.5)(2,0){5}{$\circ$}
\multiput(3.17,0.595)(2,0){4}{\line(1,0){1.86}}
\put(2.95,0.2){-$2$}
\put(4.95,0.2){-$1$}
\put(6.95,0.2){-$2$}
\put(8.95,0.2){-$1$}
\put(11,0.2){$0$}
\put(10.8,-0.2){(-$2$)}

\put(3,0.75){$1$}
\put(5,0.75){$2$}
\put(7,0.75){$5$}
\put(9,0.75){$4$}
\put(11,0.75){$3$}
\end{picture}
\vskip 0.5cm

Next, we blow up a point on the last curve, then the intersection of the last two curves, and then the intersection of the curves with $\bar{K}$ labels (-3) and (-2):

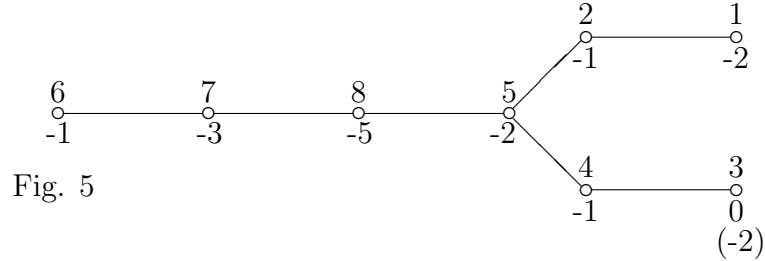
\setlength{\unitlength}{1cm}
\begin{picture}(10,2.5)
\put(0.5,-0.5){Fig. 5}
\multiput(1,0.5)(2,0){4}{$\circ$}
\multiput(1.17,0.595)(2,0){3}{\line(1,0){1.86}}
\put(0.95,0.2){-$1$}
\put(2.95,0.2){-$3$}
\put(4.95,0.2){-$5$}
\put(6.85,0.2){-$2$}

\put(1,0.75){$6$}
\put(3,0.75){$7$}
\put(5,0.75){$8$}
\put(7,0.75){$5$}

\put(7.16,0.53){\line(1,-1){0.9}}
\put(8.02,-0.53){$\circ$} 
\put(8.19,-0.42){\line(1,0){1.86}}
\put(10.02,-0.53){$\circ$}
\put(7.95,-0.84){-$1$}
\put(10.02,-0.84){$0$}
\put(9.85,-1.23){(-$2$)}
\put(8.02,-0.25){$4$}
\put(10.02,-0.25){$3$}

\put(7.16,0.67){\line(1,1){0.9}}
\put(8.02,1.51){$\circ$} 
\put(8.19,1.61){\line(1,0){1.86}}
\put(10.02,1.51){$\circ$}
\put(7.95,1.2){-$1$}
\put(9.95,1.2){-$2$}
\put(8.02,1.78){$2$}
\put(10.02,1.78){$1$}
\end{picture}
\vskip 1.8cm

Now we blow up another point on the original line at infinity, and then a point on the new curve, and again a point on the new curve, and again:

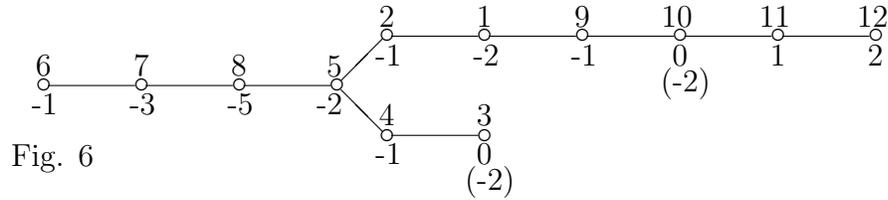
\setlength{\unitlength}{0.65cm}
\begin{picture}(10,3)
\put(0.5,-1){Fig. 6}
\multiput(1,0.5)(2,0){4}{$\circ$}
\multiput(1.25,0.64)(2,0){3}{\line(1,0){1.78}}
\put(0.92,0.05){-$1$}
\put(2.92,0.05){-$3$}
\put(4.92,0.05){-$5$}
\put(6.75,0.05){-$2$}
\put(1,0.85){$6$}
\put(3,0.85){$7$}
\put(5,0.85){$8$}
\put(6.95,0.835){$5$}

\put(7.25,0.55){\line(1,-1){0.84}}
\put(8.02,-0.53){$\circ$} 
\put(8.27,-0.38){\line(1,0){1.78}}
\put(10.02,-0.53){$\circ$}
\put(7.95,-0.97){-$1$}
\put(10.02,-0.97){$0$}
\put(9.78,-1.5){(-$2$)}
\put(8.02,-0.18){$4$}
\put(10.02,-0.18){$3$}

\put(7.25,0.75){\line(1,1){0.84}}
\multiput(8.02,1.51)(2,0){6}{$\circ$} 
\multiput(8.27,1.65)(2,0){5}{\line(1,0){1.78}}
\put(7.95,1.05){-$1$}
\put(9.95,1.05){-$2$}
\put(11.95,1.05){-$1$}
\put(14.02,1.05){$0$}
\put(16.02,1.05){$1$}
\put(18.02,1.05){$2$}
\put(13.78,0.52){(-$2$)}

\put(8.02,1.85){$2$}
\put(10.02,1.85){$1$}
\put(12.02,1.85){$9$}
\put(13.8,1.85){$10$}
\put(15.8,1.85){$11$}
\put(17.8,1.85){$12$}
\end{picture}
\vskip 1.8cm

Then, we blow up the intersection of the last two curves, and then the intersection of the curves with $\bar{K}$ labels 3 and 2:

\setlength{\unitlength}{0.65cm}
\begin{picture}(10,3)
\put(0.5,-1){Fig. 7}
\multiput(-1,0.5)(2,0){4}{$\circ$}
\multiput(-0.75,0.64)(2,0){3}{\line(1,0){1.78}}
\put(-1.08,0.05){-$1$}
\put(0.92,0.05){-$3$}
\put(2.92,0.05){-$5$}
\put(4.75,0.05){-$2$}
\put(-1,0.85){$6$}
\put(1,0.85){$7$}
\put(3,0.85){$8$}
\put(4.95,0.835){$5$}

\put(5.25,0.55){\line(1,-1){0.84}}
\put(6.02,-0.53){$\circ$} 
\put(6.27,-0.38){\line(1,0){1.78}}
\put(8.02,-0.53){$\circ$}
\put(5.95,-0.97){-$1$}
\put(8.02,-0.97){$0$}
\put(7.78,-1.5){(-$2$)}
\put(6.02,-0.18){$4$}
\put(8.02,-0.18){$3$}

\put(5.25,0.75){\line(1,1){0.84}}
\multiput(6.02,1.51)(2,0){8}{$\circ$} 
\multiput(6.27,1.65)(2,0){7}{\line(1,0){1.78}}
\put(5.95,1.05){-$1$}
\put(7.95,1.05){-$2$}
\put(9.95,1.05){-$1$}
\put(12.02,1.05){$0$}
\put(14.02,1.05){$1$}
\put(16.02,1.05){$3$}
\put(18.02,1.05){$5$}
\put(20.02,1.05){$2$}
\put(11.78,0.52){(-$2$)}

\put(6.02,1.85){$2$}
\put(8.02,1.85){$1$}
\put(10.02,1.85){$9$}
\put(11.8,1.85){$10$}
\put(13.8,1.85){$11$}
\put(15.8,1.85){$13$}
\put(17.8,1.85){$14$}
\put(19.8,1.85){$12$}

\end{picture}
\vskip 1.8cm

In what follows, we will stop keeping track of the order of creation: it is already not unique.
We will now create some branches from the curve with $\bar{K}$ label (-5) and the curve with $\bar{K}$ label (-2), adjacent to it. Specifically, we will do the following.
\begin{itemize}
\item Create 8 branches of length 1 from the (-5)-curve, by blowing up 8 distinct points on it.

\item Create 5 branches of length 3 from the (-5)-curve by blowing up a point, then the point of intersection of the new curve and the  (-5)-curve, and then another point on the (-4)-curve.

\item Create 3 forked branches from the (-5)-curve, identical to the forked branch on $Y$.

\item Create 8 branches of length 2 from the (-2)-curve, by blowing up a point, and then a point on the new curve.

\item Create 4 branches of length 3 from the (-2)-curve, by blowing up a point, then the point of intersection of the new curve and the (-2)-curve, and then another point on the (-1)-curve.

\end{itemize}

We also blow up the intersection of the curves that are indicated on the last picture by the creation numbers 2 and 5, and then the intersection of the resulting curve and the curve with the creation number 5.  Finally, we contract the strict pullback of the original line at infinity, and then contract the   curve with the creation number 9 on the last picture. The resulting surface has the following graph of curves at infinity:

\setlength{\unitlength}{0.65cm}
\begin{picture}(10,6)
\put(0.5,-1.5){Fig. 8}
\multiput(-1,0.5)(2,0){6}{$\circ$}
\multiput(-0.75,0.64)(2,0){5}{\line(1,0){1.78}}
\put(-1.08,0.05){-$1$}
\put(0.92,0.05){-$3$}
\put(2.75,0.05){-$5$}
\put(4.75,0.05){-$2$}
\put(6.92,0.05){-$1$}
\put(8.95,0.05){$0$}
\put(8.71,-0.48){(-$1$)}
\put(9.82,0.55){\circle{1}}
\put(9.4,0.36){$\times 8$}

\put(5.25,0.55){\line(1,-1){0.84}}
\put(6.02,-0.53){$\circ$} 
\put(6.27,-0.38){\line(1,0){1.78}}
\put(8.02,-0.53){$\circ$}
\put(5.95,-0.97){-$1$}
\put(8.02,-0.97){$0$}
\put(7.78,-1.5){(-$2$)}

\put(5.23,0.53){\line(1,-3){0.9}}
\multiput(6.02,-2.45)(2,0){3}{$\circ$} 
\multiput(6.27,-2.31)(2,0){2}{\line(1,0){1.78}}
\put(5.95,-2.95){-$3$}
\put(7.95,-2.95){-$1$}
\put(9.97,-2.95){$0$}
\put(9.73,-3.48){(-$1$)}
\put(10.84,-2.45){\circle{1}}
\put(10.39,-2.64){$\times 4$}

\put(5.25,0.75){\line(1,1){0.84}}
\multiput(6.02,1.51)(2,0){8}{$\circ$} 
\multiput(6.27,1.65)(2,0){7}{\line(1,0){1.78}}
\put(5.95,1.05){-$5$}
\put(7.95,1.05){-$3$}
\put(9.95,1.05){-$1$}
\put(12.02,1.05){$0$}
\put(14.02,1.05){$1$}
\put(16.02,1.05){$3$}
\put(18.02,1.05){$5$}
\put(20.02,1.05){$2$}
\put(11.78,0.52){(-$1$)}

\put(3.23,0.78){\line(1,3){0.9}}
\multiput(4.02,3.43)(2,0){7}{$\circ$}
\multiput(4.27,3.57)(2,0){6}{\line(1,0){1.78}}
\put(4.02,2.95){-$4$}
\put(5.95,2.95){-$3$}
\put(7.95,2.95){-$2$}
\put(9.95,2.95){-$1$}
\put(11.95,2.95){-$2$}
\put(13.95,2.95){-$1$}
\put(16.02,2.95){$0$}
\put(15.75,2.42){(-$2$)}
\put(12.24,3.67){\line(1,2){0.9}}
\put(13.02,5.4){$\circ$} 
\put(13.27,5.55){\line(1,0){1.78}}
\put(15.02,5.4){$\circ$}
\put(13.02,4.93){-$1$}
\put(15.02,4.93){$0$}
\put(14.78,4.4){(-$1$)}
\put(14.05,4.51){\circle{1}}
\put(13.6,4.3){$\times 3$}

\put(3.23,0.53){\line(1,-2){0.9}}
\put(4.02,-1.53){$\circ$} 
\put(3.95,-2){-$2$}
\put(3.5,-1.45){\circle{1}}
\put(3.05,-1.67){$\times 8$}

\multiput(3.15,0.75)(-1,1){3}{\line(-1,1){0.84}}
\multiput(2.06,1.51)(-1,1){3}{$\circ$}
\put(1.9,1.05){-$9$}
\put(0.9,2.05){-$4$}
\put(-0.1,3.05){-$3$}
\put(0.87,3.91){\circle{1}}
\put(0.42,3.7){$\times 5$}

\end{picture}
\vskip 1.8cm

Finally, to get the surface $Z$, we blow up some points of intersection of curves ``between" the forked (-5)-curve and (-2)-curve as on the picture below. Here the numbers above the vertices again indicate the (possible, not unique) order of creation of the new curves.

\setlength{\unitlength}{0.38cm}
\begin{picture}(10,3)
\put(-2,2.2){Fig.9}
\multiput(-1.7,0.5)(2,0){18}{$\circ$}
\multiput(-1.27,0.76)(2,0){17}{\line(1,0){1.62}}
\put(-2,-0.3){-$5$}
\put(-0.3,-0.3){-$52$}
\put(1.7,-0.3){-$47$}
\put(3.7,-0.3){-$42$}
\put(5.7,-0.3){-$37$}
\put(7.7,-0.3){-$32$}
\put(9.7,-0.3){-$27$}
\put(11.7,-0.3){-$22$}
\put(13.7,-0.3){-$39$}
\put(15.7,-0.3){-$17$}
\put(17.7,-0.3){-$12$}
\put(19.7,-0.3){-$19$}
\put(21.7,-0.3){-$26$}
\put(24,-0.3){-$7$}
\put(26,-0.3){-$9$}
\put(27.7,-0.3){-$11$}
\put(29.7,-0.3){-$13$}
\put(32,-0.3){-$2$}

\put(-0.1,1.2){10}
\put(2.3,1.2){9}
\put(4.3,1.2){8}
\put(6.3,1.2){7}
\put(8.3,1.2){6}
\put(10.3,1.2){5}
\put(12.3,1.2){4}
\put(13.9,1.2){11}
\put(16.3,1.2){3}
\put(18.3,1.2){2}
\put(19.9,1.2){12}
\put(21.9,1.2){13}
\put(24.3,1.2){1}
\put(25.9,1.2){14}
\put(27.9,1.2){15}
\put(29.9,1.2){16}
\end{picture}
\vskip 0.5cm

Clearly, the graphs for $Y$ and $Z$ describe some families of smooth compactifications of $\mathbb A^2$. Note that the curves at infinity provide a basis of the Picard group, and our labels allow us to completely describe the intersection forms. One can construct two maps between the Picard groups, $\phi_*$ and $\phi^*,$ that satisfy the projection formula and all other immediately necessary conditions for an actual Keller map $\phi$. 

Here is the general description of the map. The picture and the details for the specific branches will follow. We will generally denote the curves on $Z$ by $E_i$ and the curves on $X$ by $F_i$, where $i$ will be the $\bar{K}$ label of the curve. For all curves $E_i$ of type 2 $\phi _*(E_i)=0$; for all curves of $E_i$ of type $1$ $\phi_*(E_i)=f_iF_j$, where $f_i$ is some natural number, understood as the degree of the restriction of $\phi$ to $E_i$.  In this case the $\bar{K}$ label $i$ must be a multiple of the $\bar{K}$ label $j$: $i=e_i j$. The product $e_i\cdot f_i$ is the degree of $\phi$ in the neighborhood of the generic point on $E_i$.

{\bf General description of the map $\phi$}. The generic degree of $\phi $ is 16, so $\phi_* \circ \phi^*=16\cdot Id.$ The curve with $\bar{K}$ label 5 is of type 3, and the curve with $\bar{K}$ label 2 is of type 4. All other curves are of type 1 or 2. The forked (-5)-curve on $Z$ is sent by $\phi$ to the (-5)-curve of $Y$, with the degree of the restriction $f_{-5}=16.$ All forked  (-2)-curves on $Z$  are sent to the forked (-2)-curve on $Y$. The multi-forked (-2)-curve has $f_{-2}=13,$ and the other three (-2)-curves are mapped 1-to-1.

 The eight length 1 branches from the (-5)-curve are sent to the length 1 branch from the (-5)-curve on $Y$ with degree 2. The five length 3 branches from the (-5)-curve are sent to the length 2 branches from the (-5)-curve on $Y$ with degree 3. (Note that the degrees $e_if_i$ are constant on the branches). The  three forked branches from the (-5)-curve are sent to the forked branch from the (-5)-curve on $Y$ with degree 1. The chain between the multi-forked (-5)-curve and (-2)-curve is sent to the chain between the (-5)-curve and the (-2)-curve on $Y$, with degree 13.  The branch of length 2 from the multi-forked (-2)-curve, that ends in the curve with the $\bar{K}$ label $0$ and self-intersection (-2) is sent to the branch from the (-2)-curve on $Y$ that ends in the curve with $\bar{K}$ label $0$ and the self-intersection (-2), with degree 1. The four branches of length 3 from the multi-forked (-2)-curve are also sent there, but with degree 3. The eight branches of length 2 are sent to the branch from the (-2)-curve on $Y$ that ends in the curve with the $\bar{K}$ label $0$ and self-intersection (-1), with degree 1. Finally, the long branch that ends with the curve of type 4 is sent to the  branch from (-2)-curve on $Y$ that ends in the curve with the  $\bar{K}$ label $0$ and self-intersection (-1), with degree 5. The curve with the $\bar{K}$ label 0 on it goes to the curve with the $\bar{K}$ label 0, the curves with the $\bar{K}$ labels 1 and 3 go to some point on that curve, that the image of the curve with the $\bar{K}$ label 5 intersects.

The following picture represents the map $\phi$. To avoid overcrowding, not all arrows are drawn. 

\setlength{\unitlength}{0.65cm}
\begin{picture}(10,6)
\thicklines
\put(6,5){\bf First Framework}
\multiput(-1,0.5)(2,0){2}{\color{Plum}$\circ$}
\multiput(3,0.5)(2,0){2}{\color{red}$\circ$}
\multiput(7,0.5)(2,0){2}{\color{green}$\circ$}

\multiput(-0.75,0.64)(2,0){2}{\color{Plum}\line(1,0){1.78}}
\multiput(5.25,0.64)(2,0){2}{\color{green}\line(1,0){1.78}}
\multiput(3.26,0.62)(0.2,0){9}{\color{Cyan}.}
\put(-0.9,0.05){-$1$}
\put(0.92,0.05){-$3$}
\put(2.75,0.05){-$5$}
\put(4.6,0.05){-$2$}
\put(6.92,0.05){-$1$}
\put(8.95,0.05){$0$}
\put(8.71,-0.48){(-$1$)}
\put(9.86,0.59){\circle{1}}
\put(9.44,0.38){$\times 8$}

{\color{blue}
\put(5.25,0.55){\line(1,-1){0.84}}
\put(6.02,-0.53){$\circ$} 
\put(6.27,-0.38){\line(1,0){1.78}}
\put(8.02,-0.53){$\circ$}
}
\put(5.95,-0.97){-$1$}
\put(8.02,-0.97){$0$}
\put(7.78,-1.5){(-$2$)}

{\color{blue}
\put(5.23,0.53){\line(1,-3){0.9}}
\multiput(6.02,-2.45)(2,0){3}{$\circ$} 
\multiput(6.27,-2.31)(2,0){2}{\line(1,0){1.78}}
}
\put(5.95,-2.95){-$3$}
\put(7.95,-2.95){-$1$}
\put(9.97,-2.95){$0$}
\put(9.73,-3.48){(-$1$)}
\put(10.84,-2.45){\circle{1}}
\put(10.39,-2.64){$\times 4$}

{\color{green}
\put(5.25,0.75){\line(1,1){0.84}}
\multiput(6.02,1.51)(2,0){8}{$\circ$} 
\multiput(6.27,1.65)(2,0){7}{\line(1,0){1.78}}
}
\put(5.95,1.05){-$5$}
\put(7.95,1.05){-$3$}
\put(10,1.05){-$1$}
\put(12.02,1.05){$0$}
\put(14.02,1.05){$1$}
\put(16.02,1.05){$3$}
\put(18.02,1.05){$5$}
\put(20.02,1.05){$2$}
\put(11.78,0.52){(-$1$)}
\put(18.35,2.35){\tiny $type\  3$}
\put(18.25,1.85){$\swarrow$}

{\color{Cyan}
\put(3.23,0.78){\line(1,3){0.9}}
\multiput(4.02,3.43)(2,0){4}{$\circ$}
\multiput(4.27,3.57)(2,0){4}{\line(1,0){1.78}}
}
\put(12.02,3.43){\color{red}$\circ$}
{\color{blue}
\multiput(14.02,3.43)(2,0){2}{$\circ$}
\multiput(12.27,3.57)(2,0){2}{\line(1,0){1.78}}
}
\put(4.02,2.95){-$4$}
\put(5.95,2.95){-$3$}
\put(7.95,2.95){-$2$}
\put(9.95,2.95){-$1$}
\put(11.95,2.95){-$2$}
\put(13.95,2.95){-$1$}
\put(16.02,2.95){$0$}
\put(15.75,2.42){(-$2$)}

{\color{green}
\put(12.25,3.68){\line(1,2){0.89}}
\put(13.02,5.4){$\circ$} 
\put(13.27,5.55){\line(1,0){1.78}}
\put(15.02,5.4){$\circ$}
}
\put(13.02,4.93){-$1$}
\put(15.02,4.93){$0$}
\put(14.78,4.4){(-$1$)}
\put(14.05,4.51){\circle{1}}
\put(13.6,4.3){$\times 3$}

{\color{OliveGreen}
\put(3.23,0.53){\line(1,-2){0.89}}
\put(4.02,-1.53){$\circ$} 
}
\put(3.95,-2){-$4$}
\put(3.5,-1.45){\circle{1}}
\put(3.05,-1.67){$\times 8$}

{\color{Plum}
\multiput(3.15,0.75)(-1,1){3}{\line(-1,1){0.84}}
\multiput(2.06,1.51)(-1,1){3}{$\circ$}
}
\put(1.9,1.05){-$9$}
\put(0.9,2.05){-$4$}
\put(-0.3,3.05){-$3$}
\put(0.87,3.91){\circle{1}}
\put(0.42,3.7){$\times 5$}


\multiput(0.5,-9.52)(2,0){2}{\color{Plum} $\circ$}
\put(4.5,-9.52){\color{red} $\circ$}
\multiput(6.5,-9.52)(2,0){4}{\color{Cyan} $\circ$}
\put(14.5,-9.52){\color{red} $\circ$}
\multiput(16.5,-9.52)(2,0){2}{\color{blue}$\circ$}

\multiput(0.75,-9.37)(2,0){2}{\color{Plum} \line(1,0){1.78}}
\multiput(4.75,-9.37)(2,0){5}{\color{Cyan}\line(1,0){1.78}}
\multiput(14.75,-9.37)(2,0){2}{\color{blue}\line(1,0){1.78}}

\put(0.4,-10){-$1$}
\put(2.4,-10){-$3$}
\put(4.15,-10){-$5$}
\put(6.4,-10){-$4$}
\put(8.4,-10){-$3$}
\put(10.4,-10){-$2$}
\put(12.4,-10){-$1$}
\put(14.4,-10){-$2$}
\put(16.4,-10){-$1$}
\put(18.5,-10){$0$}
\put(18.2,-10.6){(-$2$)}

\put(4.7,-9.48){\color{OliveGreen} \line(1,-2){0.88}}
\put(5.47,-11.53){\color{OliveGreen}$\circ$} 
\put(5.4,-12){-$2$}

{\color{green}
\put(14.72,-9.28){\line(1,2){0.9}}
\put(15.5,-7.55){$\circ$} 
\put(15.75,-7.4){\line(1,0){1.78}}
\put(17.5,-7.55){$\circ$}
}
\put(15.5,-8.03){-$1$}
\put(17.5,-8.03){$0$}
\put(17.25,-8.55){(-$1$)}


{\color{Plum}
\qbezier(-0.9,0.53)(-1.5,-2)(0.42,-8.93)
\put(0.28,-9.2){$\triangleleft$}
\qbezier(0.25,3.55)(0.5,-6)(0.74,-8.78)
\put(0.5,-9.16){$\triangledown$}
}

{\color{red}
\qbezier(2.82,0.55)(0.6,-1)(4.25,-8.93)
\put(4.08,-9.2){$\triangleleft$}
}

{\color{red}
\qbezier(4.7,0.55)(5,-6)(13.78,-9.05)
\put(13.65,-9.25){\small $\triangle$}
}

{\color{OliveGreen}
\qbezier(3.64,-1.42)(5,-2)(5.04,-10.83)
\put(4.81,-11.2){$\triangledown$}
}

{\color{red}
\qbezier(11.43,3.48)(13.5,2.5)(13.8,-8.78)
\put(13.55,-9.16){$\triangledown$}
}

{\color{blue}
\qbezier(9.05,-2.38)(6,-5)(17.05,-9.05)
\put(16.95,-9.25){\small $\triangle$}
\qbezier(7.25,-0.45)(18,-4.5)(17.63,-8.8)
\put(17.4,-9.2){$\triangledown$}
\qbezier(15.24,3.51)(21,-6)(18.2,-9.05)
\put(17.8,-9.25){\small $\triangle$}
}

{\color{green}
\qbezier(7.91,0.53)(12,-4.9)(15.9,-7.05)
\put(15.75,-7.25){\small $\triangle$}
\qbezier(11.02,1.56)(15.5,-2)(16.34,-6.8)
\put(16.11,-7.2){$\triangledown$}
\qbezier(14.02,5.46)(22,-2)(16.89,-7.05)
\put(16.49,-7.25){\small $\triangle$}
}
{\color{Cyan}
\qbezier(2.7,0)(4.8,-5)(7.95,-8.95)
\put(7.8,-9.2){$\triangleleft$}
}

\put(-2,-5){\color{Plum}$deg\ 3$}
\put(1,-7){\color{red}$deg\ 16$}
\put(3.1,-8){\color{OliveGreen}$deg\ 2$}
\put(6,-8){\color{Cyan}$deg\ 13$}
\put(7,-7){\color{red}$deg\ 13$}
\put(7.5,-5){\color{blue}$deg\ 3$}
\put(13.7,-3){\color{green}$deg\ 5$}

\end{picture}
\vskip 7.8cm

\setlength{\unitlength}{0.38cm}
\begin{picture}(10,2)
\thicklines
\put(10.5,2){Close-up of the (-5)...(-2) map (light blue arrow)}
\multiput(-0.7,0.5)(34,0){12}{$\color{red} \circ$}
{\color{Cyan}
\multiput(1.3,0.5)(2,0){16}{$\circ$}
\multiput(-0.27,0.76)(2,0){17}{\line(1,0){1.62}}
}
\put(-1,-0.3){-$5$}
\put(0.6,-0.3){-$52$}
\put(2.7,-0.3){-$47$}
\put(4.7,-0.3){-$42$}
\put(6.7,-0.3){-$37$}
\put(8.7,-0.3){-$32$}
\put(10.7,-0.3){-$27$}
\put(12.7,-0.3){-$22$}
\put(14.7,-0.3){-$39$}
\put(16.7,-0.3){-$17$}
\put(18.7,-0.3){-$12$}
\put(20.7,-0.3){-$19$}
\put(22.7,-0.3){-$26$}
\put(25,-0.3){-$7$}
\put(27,-0.3){-$9$}
\put(28.7,-0.3){-$11$}
\put(30.7,-0.3){-$13$}
\put(32.7,-0.3){-$2$}

\multiput(2,-3.5)(30,0){2}{\color{red}$\circ$}
{\color{Cyan}
\multiput(8,-3.5)(6,0){4}{$\circ$}
\multiput(2.43,-3.24)(6,0){5}{\line(1,0){5.62}}
}
\put(1.8,-4.4){-$5$}
\put(7.8,-4.4){-$4$}
\put(13.8,-4.4){-$3$}
\put(19.8,-4.4){-$2$}
\put(25.8,-4.4){-$1$}
\put(31.8,-4.4){-$2$}

{\color{red}
\qbezier(-0.27,0.57)(1,-1.8)(1.8,-2.7)
\put(1.6,-3){\tiny $\triangle$}
}
{\color{Cyan}
\qbezier(1.48,0.52)(4.2,-1.5)(7.35,-2.7)
\put(7.2,-3){\tiny$\triangle$}
\qbezier(15.05,0.6)(13.65,0.2)(13.9,-2.5)
\put(13.65,-2.9){\tiny$\triangledown$}
\qbezier(23,0.6)(19.6,-1.8)(19.85,-2.5)
\put(19.6,-2.9){\tiny$\triangledown$}
\qbezier(31,0.6)(29.6,-0.5)(26.45,-2.73)
\put(25.95,-2.98){\tiny$\triangle$}
}
{\color{red}
\qbezier(32.95,0.63)(33.4,0.2)(31.95,-2.7)
\put(31.72,-3.02){\tiny$\triangleright$}
}

\put(0.5,3){Fig. 10}
\end{picture}
\vskip 2cm

We now need to fill in some details, providing $\phi_*$ and $\phi^*$ for the maps between branches of degrees more than 1. Checking the projection formula for all pairs of curves $E$ and $F$, where  $\phi(E)$ is a multiple of $F$, a point on $F,$ or a curve intersecting $F$, is left to the reader. Clearly, this would be sufficient for the projection formula, since the curves at infinity form the basis of the Picard groups of their surfaces, and for all other pairs the projection formula is trivially true.

\vskip 0.7cm
\centerline{\bf Detailed description of $\phi$}

\begin{itemize}
\item For the forked (-5)-curves on $Z$ and $Y$:

$\phi_*(E_{-5})=16F_{-5},$ $\phi^*(F_{-5})=E_{-5}$

\item For the multi-forked (-2)-curve on $Z$ and the forked (-2)-curve on $Y:$

$\phi_*(E_{-2})=13F_{-2},$ $\phi^*(F_{-2})=E_{-2}$

\item For the length 1 branches from the (-5)-curve:

$\phi_*(E_{-4})=F_{-2},$ $\phi^*(F_{-2})=2E_{-4}$

\item For the length 3 branches from the (-5)-curve:

The following picture shows where all the curves go, as well as the self-intersections of all curves, in parentheses above or below the vertices.

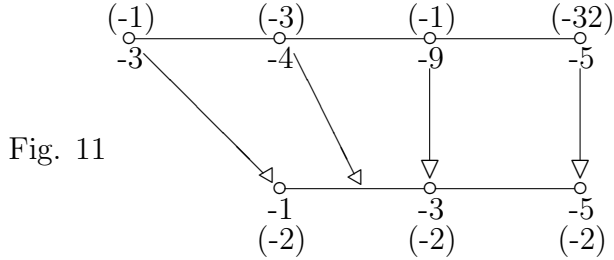
\setlength{\unitlength}{1cm}
\begin{picture}(10,1.5)
\put(-0.5,-1){Fig. 11}
\multiput(1,0.5)(2,0){4}{$\circ$}
\multiput(1.17,0.595)(2,0){3}{\line(1,0){1.86}}
\multiput(3,-1.5)(2,0){3}{$\circ$}
\multiput(3.17,-1.4)(2,0){2}{\line(1,0){1.86}}

\put(0.95,0.2){-$3$}
\put(2.95,0.2){-$4$}
\put(4.95,0.2){-$9$}
\put(6.95,0.2){-$5$}
\put(0.8,0.75){(-$1$)}
\put(2.8,0.75){(-$3$)}
\put(4.8,0.75){(-$1$)}
\put(6.7,0.75){(-$32$)}

\put(2.95,-1.8){-$1$}
\put(4.95,-1.8){-$3$}
\put(6.95,-1.8){-$5$}
\put(2.8,-2.2){(-$2$)}
\put(4.8,-2.2){(-$2$)}
\put(6.8,-2.2){(-$2$)}

\put(4,-1.35){$\triangleleft$}
\put(2.82,-1.32){$\triangleleft$}
\put(4.95,-1.26){$\triangledown$}
\put(6.95,-1.26){$\triangledown$}

\put(1.3,0.4){\line(1,-1){1.59}}
\put(3.3,0.4){\line(1,-2){0.8}}
\put(5.1,0.2){\line(0,-1){1.23}}
\put(7.1,0.2){\line(0,-1){1.23}}

\end{picture}
\vskip 2.3cm

$\phi_*(E_{-9})=F_{-3},$ $\phi_*(E_{-4})=0,$ $\phi_*(E_{-3})=F_{-1},$ 

$\phi^*(F_{-3})=3E_{-9}+ E_{-4},$ $\phi^*(F_{-1})=3E_{-3}+ E_{-4}$ 

\item  For the length 3 branches from the (-2)-curve:

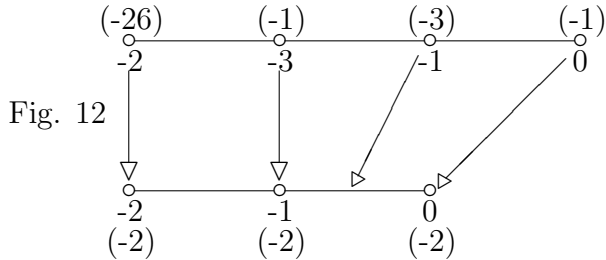
\setlength{\unitlength}{1cm}
\begin{picture}(10,1.5)
\put(-0.5,-0.5){Fig. 12}
\multiput(1,0.5)(2,0){4}{$\circ$}
\multiput(1.17,0.595)(2,0){3}{\line(1,0){1.86}}
\multiput(1,-1.5)(2,0){3}{$\circ$}
\multiput(1.17,-1.4)(2,0){2}{\line(1,0){1.86}}

\put(0.95,0.2){-$2$}
\put(2.95,0.2){-$3$}
\put(4.95,0.2){-$1$}
\put(7,0.2){$0$}
\put(0.7,0.75){(-$26$)}
\put(2.8,0.75){(-$1$)}
\put(4.8,0.75){(-$3$)}
\put(6.8,0.75){(-$1$)}

\put(0.95,-1.8){-$2$}
\put(2.95,-1.8){-$1$}
\put(5,-1.8){$0$}
\put(0.8,-2.2){(-$2$)}
\put(2.8,-2.2){(-$2$)}
\put(4.8,-2.2){(-$2$)}

\put(0.95,-1.26){$\triangledown$}
\put(2.95,-1.26){$\triangledown$}
\put(4.05,-1.35){$\triangleright$}
\put(5.2,-1.37){$\triangleright$}

\put(1.1,0.2){\line(0,-1){1.23}}
\put(3.1,0.2){\line(0,-1){1.23}}
\put(4.95,0.4){\line(-1,-2){0.8}}
\put(6.91,0.36){\line(-1,-1){1.59}}

\end{picture}
\vskip 2.3cm

$\phi_*(E_{-3})=F_{-1},$ $\phi_*(E_{-1})=0,$ $\phi_*(E_{0})=F_{0},$ 

$\phi^*(F_{-1})=3E_{-3}+ E_{-1},$ $\phi^*(F_{0})=3E_{0}+ E_{-1}$ 

\item For the long branch from the (-2)-curve:

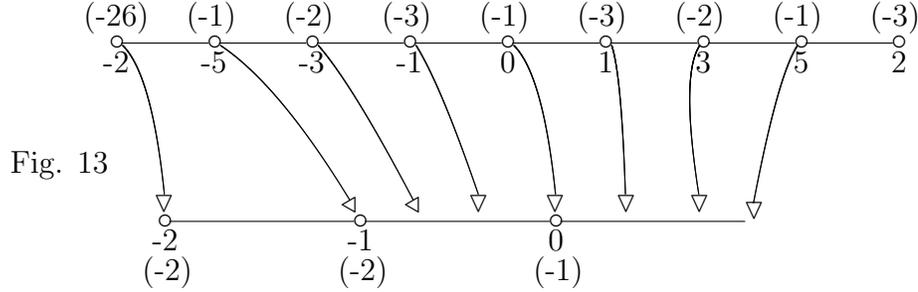
\setlength{\unitlength}{0.65cm}
\begin{picture}(10,3.5)
\put(-3,-1){Fig. 13}
\multiput(-0.98,1.51)(2,0){9}{$\circ$} 
\multiput(-0.73,1.65)(2,0){8}{\line(1,0){1.78}}
\put(-1.1,1.05){-$2$}
\put(0.9,1.05){-$5$}
\put(2.9,1.05){-$3$}
\put(4.9,1.05){-$1$}
\put(7.02,1.05){$0$}
\put(9.02,1.05){$1$}
\put(11.02,1.05){$3$}
\put(13.02,1.05){$5$}
\put(15.02,1.05){$2$}

\put(-1.5,2){(-$26$)}
\put(0.6,2){(-$1$)}
\put(2.6,2){(-$2$)}
\put(4.6,2){(-$3$)}
\put(6.6,2){(-$1$)}
\put(8.6,2){(-$3$)}
\put(10.6,2){(-$2$)}
\put(12.6,2){(-$1$)}
\put(14.6,2){(-$3$)}

\multiput(0,-2.14)(4,0){3}{$\circ$} 
\multiput(0.25,-2)(4,0){3}{\line(1,0){3.78}}
\put(-0.1,-2.6){-$2$}
\put(3.9,-2.6){-$1$}
\put(8,-2.6){$0$}
\put(-0.3,-3.2){(-$2$)}
\put(3.7,-3.2){(-$2$)}
\put(7.7,-3.2){(-$1$)}

\put(-0.07,-1.8){\small $\triangledown$}
\put(3.76,-1.81){$\triangleleft$}
\put(5.06,-1.81){$\triangleleft$}
\put(6.36,-1.8){\small $\triangledown$}
\put(7.93,-1.8){\small $\triangledown$}
\put(9.38,-1.8){\small $\triangledown$}
\put(10.88,-1.8){\small $\triangledown$}
\put(12,-1.95){\small $\triangledown$}

\qbezier(-0.7,1.6)(-0.1,1)(0.15,-1.44)
\qbezier(1.3,1.6)(2.5,0.8)(3.94,-1.56)
\qbezier(3.3,1.6)(4,0.8)(5.24,-1.56)
\qbezier(5.3,1.6)(5.8,0.8)(6.57,-1.44)
\qbezier(7.3,1.6)(7.9,1)(8.15,-1.44)
\qbezier(9.3,1.6)(9.5,1)(9.58,-1.44)
\qbezier(11.1,1.6)(10.7,1)(11.08,-1.44)
\qbezier(13.1,1.6)(12.7,1)(12.2,-1.62)

\end{picture}
\vskip 3cm
$\phi_*(E_{-5})=F_{-1},$   $\phi_*(E_{-3})=\phi_*(E_{-1})=0,$ $\phi_*(E_{0})=F_{0},$

$\phi_*(E_{1})=\phi_*(E_{3})=0,$  $\phi (E_5)$ is a curve, generically in $\mathbb A^2,$

$\phi(E_2)$ is a point in $\mathbb A^2,$

$\phi^*(F_{-1})=5E_{-5}+3E_{-3}+ E_{-1},$

$\phi^*(F_{0})=5E_{0}+ (2E_{-1}+E_{-3})+(2E_1+E_3)$ 

\item Finally, for the chain between the (-5)-curve and the (-2)-curve:

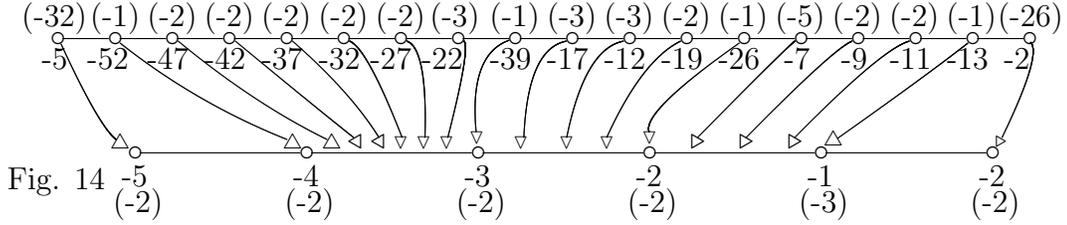
\setlength{\unitlength}{0.38cm}
\hskip -1.3cm \begin{picture}(10,3)
\put(-2.2,-4.5){Fig. 14}
\multiput(-0.7,0.5)(34,0){12}{$\circ$}

\multiput(1.3,0.5)(2,0){16}{$\circ$}
\multiput(-0.27,0.76)(2,0){17}{\line(1,0){1.62}}

\put(-1,-0.3){-$5$}
\put(0.6,-0.3){-$52$}
\put(2.7,-0.3){-$47$}
\put(4.7,-0.3){-$42$}
\put(6.7,-0.3){-$37$}
\put(8.7,-0.3){-$32$}
\put(10.5,-0.3){-$27$}
\put(12.3,-0.3){-$22$}
\put(14.7,-0.3){-$39$}
\put(16.7,-0.3){-$17$}
\put(18.7,-0.3){-$12$}
\put(20.7,-0.3){-$19$}
\put(22.7,-0.3){-$26$}
\put(25,-0.3){-$7$}
\put(27,-0.3){-$9$}
\put(28.7,-0.3){-$11$}
\put(30.7,-0.3){-$13$}
\put(32.7,-0.3){-$2$}

\put(-1.7,1.2){(-$32$)}
\put(0.7,1.2){(-$1$)}
\put(2.7,1.2){(-$2$)}
\put(4.7,1.2){(-$2$)}
\put(6.7,1.2){(-$2$)}
\put(8.7,1.2){(-$2$)}
\put(10.7,1.2){(-$2$)}
\put(12.6,1.2){(-$3$)}
\put(14.7,1.2){(-$1$)}
\put(16.7,1.2){(-$3$)}
\put(18.7,1.2){(-$3$)}
\put(20.7,1.2){(-$2$)}
\put(22.7,1.2){(-$1$)}
\put(24.7,1.2){(-$5$)}
\put(26.7,1.2){(-$2$)}
\put(28.7,1.2){(-$2$)}
\put(30.7,1.2){(-$1$)}
\put(32.5,1.2){(-$26$)}

\multiput(2,-3.5)(30,0){2}{$\circ$}
\multiput(8,-3.5)(6,0){4}{$\circ$}
\multiput(2.43,-3.24)(6,0){5}{\line(1,0){5.62}}

\put(1.8,-4.4){-$5$}
\put(7.8,-4.4){-$4$}
\put(13.8,-4.4){-$3$}
\put(19.8,-4.4){-$2$}
\put(25.8,-4.4){-$1$}
\put(31.8,-4.4){-$2$}

\put(1.5,-5.3){(-$2$)}
\put(7.5,-5.3){(-$2$)}
\put(13.5,-5.3){(-$2$)}
\put(19.5,-5.3){(-$2$)}
\put(25.5,-5.3){(-$3$)}
\put(31.5,-5.3){(-$2$)}

\qbezier(-0.32,0.57)(1,-2.3)(1.6,-2.77)
\put(1.4,-3.1){\tiny $\triangle$}
\qbezier(1.75,0.6)(4.4,-1.5)(7.58,-2.77)
\put(7.43,-3.1){\tiny$\triangle$}
\qbezier(15.35,0.7)(13.95,0.2)(14.2,-2.5)
\put(13.95,-2.9){\tiny$\triangledown$}
\qbezier(23.4,0.6)(20,-1.8)(20.25,-2.5)
\put(20,-2.9){\tiny$\triangledown$}
\qbezier(31.4,0.6)(30,-0.5)(26.85,-2.73)
\put(26.35,-2.98){\tiny$\triangle$}
\qbezier(33.6,0.63)(34,0.2)(32.6,-2.7)
\put(32.37,-3.02){\tiny$\triangleright$}

\qbezier(3.78,0.67)(6.4,-1.5)(8.98,-2.77)
\put(8.8,-3.1){\tiny$\triangle$}
\qbezier(5.78,0.67)(8.4,-1.5)(9.91,-2.7)
\put(9.7,-3.1){\small$\triangleleft$}

\qbezier(7.78,0.67)(9.8,-1.1)(10.8,-2.7)
\put(10.5,-3.1){\small$\triangleleft$}
\qbezier(9.78,0.67)(10.4,0.4)(11.61,-2.7)
\put(11.3,-3.1){\tiny $\triangledown$}
\qbezier(11.78,0.67)(12.5,0.4)(12.41,-2.7)
\put(12.1,-3.1){\tiny $\triangledown$}
\qbezier(13.78,0.67)(14,0.4)(13.21,-2.7)
\put(12.9,-3.1){\tiny $\triangledown$}

\qbezier(17.35,0.67)(16,0.4)(15.87,-2.7)
\put(15.5,-3.1){\tiny $\triangledown$}
\qbezier(19.35,0.67)(18,0.4)(17.37,-2.7)
\put(17.1,-3.1){\tiny $\triangledown$}
\qbezier(21.35,0.67)(20,0)(18.77,-2.7)
\put(18.5,-3.1){\tiny $\triangledown$}

\qbezier(25.35,0.67)(25,0.4)(21.95,-2.7)
\put(21.7,-3.1){\small$\triangleright$}
\qbezier(27.35,0.67)(26.7,0.4)(23.65,-2.7)
\put(23.4,-3.1){\small$\triangleright$}
\qbezier(29.35,0.67)(28.4,0.4)(25.35,-2.7)
\put(25.1,-3.1){\small$\triangleright$}
\end{picture}
\vskip 2.3cm
$\phi_*(E_{-52})=F_{-4},$   $\phi_*(E_{-47})=\phi_*(E_{-42})=\phi_*(E_{-37})=$

$\phi_*(E_{-32})=\phi_*(E_{-27})=\phi_*(E_{-22})=0,$ $\phi_*(E_{-39})=F_{-3},$

$\phi_*(E_{-17})=\phi_*(E_{-12})=\phi_*(E_{-19})=0,$ $\phi_*(E_{-26})=F_{-2},$

$\phi_*(E_{-7})=\phi_*(E_{-9})=\phi_*(E_{-11})=0,$ $\phi_*(E_{-13})=F_{-1}$,

\hskip -0.8cm $\phi^*(F_{-4})=13E_{-52}+11E_{-47}+ 9E_{-42}+7E_{-37}+5E_{-32}+3E_{-27}+E_{-22},$

\hskip -0.8cm $\phi^*(F_{-3})=13E_{-39}+(E_{-47}+ 2E_{-42}+3E_{-37}+4E_{-32}+5E_{-27}+6E_{-22})+
(5E_{-17}+2E_{-12}+E_{-19}),$

\hskip -0.8cm $\phi^*(F_{-2})=13E_{-26}+(E_{-17}+3E_{-12}+8E_{-19})+(3E_{-7}+2E_{-9}+E_{-11}),$

\hskip -0.8cm $\phi^*(F_{-1})=13E_{-13}+9E_{-11}+5E_{-9}+E_{-7}$
\end{itemize}

\section{First Framework, Continued}

We will now dig deeper into our framework. For the rest of this section, we will assume that it corresponds to the actual map $\phi$ (even though it probably does not).

First, you may have noticed that we have two Belyi maps: the red arrows in the main picture, above the forked (-5)-curve and (-2)-curve.  One can draw the corresponding rational dessin d'enfants and, furthermore, write down the explicit rational functions that define them. (Thanks to Maple for their Gr\"obner basis implementation and more!) Then we will find  the degrees of the polynomials that define $\phi$. Finally, we will write down explicitly the coordinates on some open subsets of $Z$ and $Y$  and explain how to use them to search for $\phi$.

{\bf Belyi map of (-5)-curves, the First Framework}

The degree of the map is 16. The map is ramified above three points: intersection with the (-2)-curve, (-3)-curve, and (-4)-curve, that we will identify with $\{0\}$, $\{\infty\},$ and $\{1\}$ respectively. Above $\{0\}$, we have 8 ramification points of order 2, so the corresponding dessin is ``clean". Above $\{\infty\}$ we have 5 points of order 3 and 1 point of order 1, and above $\{1\}$ we have one point of order 13 and 3 points of order 1. A simple combinatorial analysis leads to the following clean dessin d'enfant (unique as a graph, but not as a dessin):

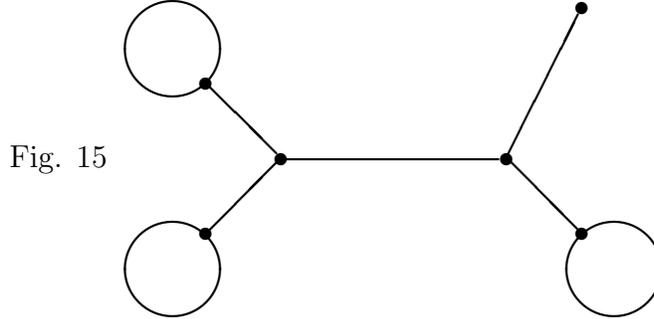
\setlength{\unitlength}{1cm}
\begin{picture}(10,2.2)
\put(0.5,0){Fig. 15}
\thicklines
\multiput(4,0)(3,0){2}{$\bullet$}
\put(3,1){$\bullet$}
\put(8,2){$\bullet$}
\multiput(3,-1)(5,0){2}{$\bullet$}
\put(4.15,0.1){\line(3,0){2.86}}
\put(3.15,-0.85){\line(1,1){0.9}}
\put(7.13,0.17){\line(1,2){0.93}}
\multiput(3.15,1.08)(4,-1){2}{\line(1,-1){0.9}}
\put(2.67,1.57){\circle{1.3}}
\put(8.54,-1.36){\circle{1.3}}
\put(2.67,-1.36){\circle{1.3}}
\end{picture}
\vskip 1.9cm

If we parametrize the (-5)-curve on $Z$ by some parameter $w$ that equals $0$ at the unique point of order 1 above $\{\infty\}$ and equals $\infty$ at the point of order $13$ above $\{1\}$, we will only have one degree of freedom left: multiplying $w$ by a non-zero number. Up to that, our Belyi map is given (thanks to Maple) as  $w\mapsto \frac{p^2(w)}{w\cdot r^3(w)},$ where the polynomials $p(w)$ and $r(w)$ can be chosen as follows:
$$p(w)= w^8+(2+8\sqrt{-3})w^7+\frac{-233+50\sqrt{-3}}{3}w^6+\frac{-4600-376\sqrt{-3}}{3}w^5+$$
$$\frac{835-890\sqrt{-3}}{3}w^4+\frac{2420+22\sqrt{-3}}{3}w^3+(\frac{1043}{3}+336\sqrt{-3})w^2+$$
$$(-118+158\sqrt{-3})w+(-28+41\sqrt{-3}),$$

$$r(w)=w^5+\frac{4+16\sqrt{-3}}{3}w^4+\frac{-278+68\sqrt{-3}}{9}w^3+(-\frac{140}{3}-24\sqrt{-3})w^2+$$
$$\frac{35-112\sqrt{-3}}{3}w+\frac{68-20\sqrt{-3}}{3}.$$

Note that $deg(p^2-w\cdot r^3)=3$.

{\bf Belyi map of (-2)-curves, the First Framework}

The degree of the map is 13. The map is ramified above three points. They correspond to the branches that end with the 0-curve with the self-intersection  (-2), with the 0-curve with the self-intersection (-1), and with the forked (-5)-curve. We will identify them as $\{0\}$,  $\{1\}$, and $\{\infty\}$ respectively. On the (-2)-curve on $Z$ we will identify with $t=\infty$ the unique point that is sent to $\{\infty\}$, so that the Belyi map is given by a polynomial. Then above $\{0\}$ we have one point of order 1, that we will call $t=0$, and 4 points of order 3. Above $\{1\}$ we have one point of order 5, that we will call $t=1,$ and 8 points of order 1. The corresponding unique dessin d'enfant is the following:

\setlength{\unitlength}{1cm}
\begin{picture}(10,1.7)
\put(0.5,-1){Fig. 16}
\thicklines
\multiput(3,0)(3,0){3}{$\bullet$}
\multiput(4,0)(4,0){2}{$\circ$}
\multiput(3.17,0.1)(5,0){2}{\line(1,0){0.86}}
\multiput(4.17,0.1)(2,0){2}{\line(1,0){1.86}}
\put(6.11,0.15){\line(0,1){1.18}}
\put(6,1.28){$\circ$}
\put(6.11,0.1){\line(1,-2){0.96}}
\put(6.11,0.1){\line(-1,-2){0.96}}
\put(5,-2){$\circ$}
\put(7,-2){$\circ$}
\multiput(4,-2)(4,0){2}{$\bullet$}
\multiput(4.17,-1.9)(3,0){2}{\line(1,0){0.86}}
\multiput(4.05,0.15)(1,-2){2}{\line(-1,1){0.6}}
\multiput(8.16,0.15)(-1,-2){2}{\line(1,1){0.6}}
\multiput(3.3,0.7)(1,-2){2}{$\bullet$}
\multiput(8.7,0.7)(-1,-2){2}{$\bullet$}
\end{picture}
\vskip 2.5cm
The Belyi map, for our choice of $t$, is given by 
$$t\mapsto \frac{1}{3^{15}}t(35t^4-182t^3+390t^2-455t+455)^3$$

In what follows, we will use $(x_1,x_2)$ as the pair of coordinates on the source $\mathbb A^2$ and $(y_1,y_2)$ as the coordinates on the target $\mathbb A^2.$ So $x_1$ and $x_2$ are rational functions on $Z$, while $y_1$ and $y_2$ are rational functions on $Y$. Linear transformations of the source and target planes give us some freedom on what points to blow up, so our first blowup on $Y$ will be the point of intersection of the line at infinity and the line $y_2=0.$ On the (-1)-curve created by that blowup  $y_1$ has a pole of order $1$ and $y_2$ is a parameter, its valuation is $0$. On the next curve, with $\bar{K}$ label (-3), the valuations of $y_1$ and $y_2$ are (-2) and (-1); on the (-5)-curve they are (-3) and (-2).  So $\frac{y_1^2}{y_2^3}$ is a parameter on the (-5)-curve. By scaling the coordinates, we may assume that the next blowup in the creation of $Y$ was at the point $\frac{y_1^2}{y_2^3}=1.$ Valuations of $y_1$ and $y_2$ with respect to this (-4)-curve are also (-3) and (-2). The same will be true all the way to the curve with the $\bar{K}$ valuation 0. Then for the (-1)-curve next to it the valuations are (-6) and (-4), and for the forked (-2)-curve and the remaining two curves on $Y$ they are (-9) and (-6).

From our description of $\phi$ we can get the valuations of $\phi^*(y_1)$ and  $\phi^*(y_2)$. Specifically, from the formulas under Figure 13, the valuations of $\phi^*(y_1)$ and  $\phi^*(y_2)$ with the respect to the 0-curve are (-45) and (-30). The valuations of $\phi^*(y_1)$ and  $\phi^*(y_2)$ with respect to the (-1)-curve are (-27) and (-18) (note that $y_1$ and $y_2$ have poles at $F_{-1}$ and $F_0,$ and they both contribute). Now recall that to get back to the original line at infinity on $\mathbb P^2$ we need to blow up twice. First, we get a curve with $\bar{K}$ label (-1), for which the valuations are (-72) and (-48), and then the original line at infinity on $\mathbb P^2$ with valuations (-99) and (-66). If we choose the coordinates $(x_1,x_2)$ on the source $\mathbb A^2$ judicially, $y_1$ will be a polynomial in $x_1$ and $x_2$ of separate degrees 27 and 72 and total degree 99. And $y_2$ will be a polynomial in $x_1$ and $x_2$ of separate degrees 18 and 48 and total degree 66.

\begin{Remark} Back in 1983 T.-T. Moh published a proof that there are no Keller maps of degrees up to 100 (cf. \cite{Moh}). In fact, the (99,66) pair of degrees was the last troublesome case that he had to discard. Unfortunately, Moh's paper does not give full details of his proof for that case. Only a reduction to an easier problem is presented, and the argument is sketchy: it seems that in order to understand it, one has to fully understand the argument for some smaller pairs of degrees, - something that I was unable to do.  To complicate the matter further,  in 2016 Yansong Xu, a 1993 Ph.D. student of Moh, posted a preprint in which he claimed that Moh's proof had a gap, that he just managed to patch. According to Christian Valqui, Yansong Xu's argument had a mistake, that he acknowledged (cf. \cite{GGHV}), but Rodrigo Horruitiner essentially proved it in his Master's thesis. Finally, my own calculations, using Maple, based on the ideas below, led to the same result: no map. So, in all likelihood, there is no map $\phi$ that satisfies our framework, but we currently do not have a {\bf simple} reason for this.
\end{Remark}

In order to fully understand the structure of $Z$ (and $Y$), we need some notation for the local coordinates of some Zariski open subsets.

\begin{Definition} Suppose $E_i$ and  $E_j$ are two intersecting curves at infinity. This means that on the graph of curves, we have an edge connecting $E_i$ and $E_j$. We will call a pair of rational functions on our surface $(f_i,f_j)$ a local coordinate system for that edge if
$$\left\{
\begin{array}{l}
(f_i)_{|_{E_i}}:E_i \to \mathbb P^1\  is\  an\  isomorphism\\
(f_j)_{|_{E_j}}:E_j \to \mathbb P^1\  is\  an\  isomorphism\\
(f_i)\ has\ zero\ on\ {E_j}\ (of\ order\ 1)\\
(f_j)\ has\ zero\ on\ {E_i}\ (of\ order\ 1)
\end{array}
\right.$$
\end{Definition}

Note that the coordinate system of an edge is far from unique, and that $f_i$ may have other zeros and poles, possibly intersecting $E_j$, and the same for $f_j$ and $E_i.$ On the other hand, a local coordinate system exists for every edge of our graph of curves, which can be proven by induction. We will keep track of this notationally as follows:

\setlength{\unitlength}{1cm}
\begin{picture}(10,1)
\put(0.5,0){Fig. 17}
\multiput(5,0)(2,0){2}{$\circ$}
\put(5.15,0.13){\line(1,0){1.86}}
\put(5.1,0.33){$(f_i\ \ \  ,\ \ \ f_j)$}
\put(4.9,-0.35){$E_i$}
\put(6.9,-0.35){$E_j$}
\end{picture}
\vskip 1cm

When an edge $E_iE_j$ is ``broken" by the blowup of $E_i\cap E_j,$ to create $E_k,$ we get the following:

\setlength{\unitlength}{1cm}
\begin{picture}(10,1)
\put(0.5,0){Fig. 18}
\multiput(4,0)(2,0){3}{$\circ$}
\multiput(4.15,0.13)(2,0){2}{\line(1,0){1.86}}
\put(4.1,0.37){$(f_i\ \ \  ,\ \ \ \frac{f_j}{f_i})$}
\put(6.1,0.37){$(\frac{f_i}{f_j}\ \ \  ,\ \ \ f_j)$}
\put(3.9,-0.35){$E_i$}
\put(5.9,-0.35){$E_k$}
\put(7.9,-0.35){$E_j$}
\end{picture}
\vskip 1cm

When a non-intersection point is blown up on a curve at infinity, the procedure is more complicated. Since at this time we must have some parameter on the curve (a rational function, identifying it with $\mathbb P^1$), subtracting a constant from that parameter will give us a function that has a zero of order one on some curve intersecting our curve transversally at any given point. We also need, however, a function that equals zero on our curve, and {\it has no other zeros or poles passing through the point to be blown up}. For our purposes this can be done by hand, but it would be nice to automatize this. Non-uniqueness of the pair of coordinates is part of the problem here.

As an illustration, below are (some possible) coordinate systems for the surface in Figure 7.  Here $e_1$ and $e_2$ are arbitrary parameters.

\setlength{\unitlength}{0.65cm}
\begin{picture}(10,3.8)
\put(0,2.5){Fig. 19}
\multiput(-1,0.5)(2,0){4}{$\circ$}
\multiput(-0.75,0.64)(2,0){3}{\line(1,0){1.78}}
\put(-1.08,0.05){-$1$}
\put(0.92,0.05){-$3$}
\put(2.92,0.05){-$5$}
\put(4.75,0.05){-$2$}

\put(5.25,0.55){\line(1,-1){0.84}}
\put(6.02,-0.53){$\circ$} 
\put(6.27,-0.38){\line(1,0){2.48}}
\put(8.72,-0.53){$\circ$}
\put(5.95,-0.97){-$1$}
\put(8.72,-0.97){$0$}
\put(8.38,-1.5){(-$2$)}

\put(5.25,0.75){\line(1,1){0.84}}
\multiput(6.02,1.51)(2,0){8}{$\circ$} 
\multiput(6.27,1.65)(2,0){7}{\line(1,0){1.78}}
\put(5.95,1.05){-$1$}
\put(7.95,1.05){-$2$}
\put(9.95,1.05){-$1$}
\put(12.02,1.05){$0$}
\put(14.02,1.05){$1$}
\put(16.02,1.05){$3$}
\put(18.02,1.05){$5$}
\put(20.02,1.05){$2$}


\put(6.2,2){$(\frac{1}{x_2},\frac{x_2}{x_1})$}
\put(8.2,2){$(\frac{x_1}{x_2},\!\frac{1}{x_1})$}
\put(10.15,2){$(\!x_1,\!\!\frac{1}{x_1x_2}\!)$}
\put(3.8,2.3){$(\frac{1}{x_1x_2^3},x_2)$}
\put(5.4,1.5){$\downarrow$}

\put(6.1,0){$(x_1x_2^2,\!\frac{1}{x_1x_2})$}
\put(3.8,-1.5){$(x_1x_2^3,\frac{1}{x_1x_2^2})$}
\put(5.4,-0.3){$\uparrow$}
\qbezier(5.565,-1)(5.565,-0.8)(5.565,-0.3)

\put(-0.5,1.5){$(\frac{x_2}{(x_1x_2^3-1)^3},x_1x_2^3-1)$}
\put(3.5,0.8){$\searrow$}

\put(-0.8,-0.7){$(\!\frac{x_2}{(x_1x_2^3-1)^2},\! \!\frac{1}{x_2}(x_1x_2^3-1)^3\!)$}
\put(2,0){$\uparrow$}

\put(-1,-2.3){$(\!\frac{x_2}{x_1x_2^3-1},\! \!\frac{1}{x_2}(x_1x_2^3-1)^2\!)$}
\put(-0.7,0){$\nearrow$}
\qbezier(-0.3,-1.8)(-2,-1.3)(-0.64,-0.1)

\put(9.2,0.5){$(\!x_1x_2\!-\!e_1,\!\!\frac{1}{x_1x_2^2-e_1x_2}\!)$}
\put(12.5,1.1){$\nearrow$}

\put(9.2,2.8){$(x_1x_2^2\!-\!e_1x_2\!-\!e_2,\frac{x_2}{(x_1x_2^3-e_1x_2^2-e_2x_2)^2})$}
\put(14.5,1.9){$\searrow$}

\put(9.5,-0.7){$(\!\frac{(x_1x_2^3-e_1x_2^2-e_2x_2)^2}{x_2},\!\frac{x_2}{(x_1x_2^3-e_1x_2^2-e_2x_2)^3}\!)$}
\put(16.5,1.1){$\nearrow$}
\qbezier(15.2,-0.2)(15.7,0.2)(16.55,1)

\put(11.5,-1.9){$(\!\frac{(x_1x_2^3-e_1x_2^2-e_2x_2)^3}{x_2},\!\frac{1}{(x_1x_2^3-e_1x_2^2-e_2x_2)}\!)$}
\put(19,1.1){$\uparrow$}
\qbezier(18.2,-1.2)(19,-0.3)(19.167,1.05)

\end{picture}
\vskip 1.8cm

These edge coordinate systems can be used to search for a map $\phi$ as follows. Suppose $y_1=\sum \limits_{i=0}^{27} \sum \limits_{j=0}^{72} a_{ij}x_1^ix_2^j$, $y_2=\sum \limits_{i=0}^{18} \sum \limits_{j=0}^{48} b_{ij}x_1^ix_2^j.$ Take, for example, an edge between the curves with $\bar{K}$ labels 0 and 1:
$$(\alpha, \beta)=(x_1x_2-e_1,\frac{1}{x_1x_2^2-e_1x_2})$$
We can solve for $(x_1,x_2)$ in terms of $(\alpha,\beta):$
$$(x_1,x_2)=( \alpha^2\beta+e_1\alpha \beta, \frac{1}{\alpha \beta})$$
Plugging that into the formulas for $y_1$ and $y_2$ and using the known valuations of $\phi^*y_1$ and $\phi^*y_2,$ we get many equations on coefficients $a_{ij}$ and $b_{ij}$. These equations are linear in the coefficients, but non-linear in $e_1$ and $e_2.$ Altogether, there are hundreds of variables and hundreds of equations. In actuality, I was writing $y_1$ and $y_2$ as Laurent polynomials in $v=x_1x_2^3-1$ and $w=\frac{1}{x_2}(x_1x_2^3-1)^3,$ that are parameters on the (-2)-curve and (-5)-curve respectively. After solving, using Maple, the hundreds of linear equations on hundreds of variables, I got down to just a dozen or so coefficients, and with a bit more work figured out that no map $\phi $ can exist. It should be stressed that one careless mistake anywhere in the process would likely lead to a missed solution, and I cannot trust my own bookkeeping abilities to claim that I actually have a proof that no $\phi $ exists.

Not that we are too far off target: here is a rather interesting map. If  $p$ and $r$ are the two polynomials from the degree 16 Belyi map, consider $(x_1,x_2)\mapsto (y_1,y_2)$, where
$$ \left\{
\begin{array}{l}
y_1=x_1^3x_2^8\cdot p(\frac{1}{x_2}(x_1x_2^3-1)^3)\\
\ \\
y_2=x_1^2x_2^5(x_1x_2^3-1)\cdot r(\frac{1}{x_2}(x_1x_2^3-1)^3)
\end{array}\right.$$

Note that this map is polynomial, and the polynomials have the correct degrees. It has a rather simple Jacobian, a constant multiple of $x_1^4x_2^{12}$. It is generically 16:1. It has some other required properties for a Keller map, for example it is not proper (there are curves of type 3 and 4, in fact exactly those that our framework predicts, for $e_1=e_2=0$). Unfortunately, there does not seem to be a way to modify it to get a Keller map.

\section{Second Framework}

Since we seem to be out of luck with the above framework, it makes sense to look for more sophisticated ones. Indeed, one can be obtained if when constructing $Y$ we create one more (-1)-curve from  the 0-curve. Specifically, here is a graph for $Y$:

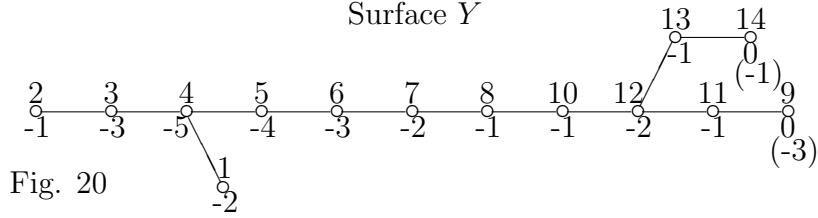
\setlength{\unitlength}{.5cm}
\begin{picture}(10,4)
\put(10,3){Surface $Y$}
\put(1,-1.5){Fig. 20}
\multiput(1.5,0.48)(2,0){11}{$\circ$} 
\multiput(1.82,0.67)(2,0){10}{\line(1,0){1.71}}
\put(1.4,0){-$1$}
\put(3.4,0){-$3$}
\put(5.1,0){-$5$}
\put(7.4,0){-$4$}
\put(9.4,0){-$3$}
\put(11.4,0){-$2$}
\put(13.4,0){-$1$}
\put(15.4,0){-$1$}
\put(17.4,0){-$2$}
\put(19.4,0){-$1$}
\put(21.5,0){$0$}
\put(21.2,-0.6){(-$3$)}

\put(1.5,0.9){$2$}
\put(3.5,0.9){$3$}
\put(5.5,0.9){$4$}
\put(7.5,0.9){$5$}
\put(9.5,0.9){$6$}
\put(11.5,0.9){$7$}
\put(13.5,0.9){$8$}
\put(15.3,0.9){$10$}
\put(17.05,0.9){$12$}
\put(19.3,0.9){$11$}
\put(21.5,0.9){$9$}

\put(5.76,0.53){\line(1,-2){0.86}}
\put(6.47,-1.53){$\circ$} 
\put(6.4,-2){-$2$}
\put(6.51,-1.1){$1$}

\put(17.77,0.82){\line(1,2){0.86}}
\put(18.5,2.45){$\circ$} 
\put(18.83,2.65){\line(1,0){1.71}}
\put(20.5,2.45){$\circ$}
\put(18.3,2.87){$13$}
\put(20.3,2.87){$14$}
\put(18.5,1.97){-$1$}
\put(20.5,1.97){$0$}
\put(20.25,1.45){(-$1$)}
\end{picture}
\vskip 1.8cm

The surface $Z$ is constructed similarly to the surface $Z$ in the first framework. Specifically, here is a new version of Figure 7:

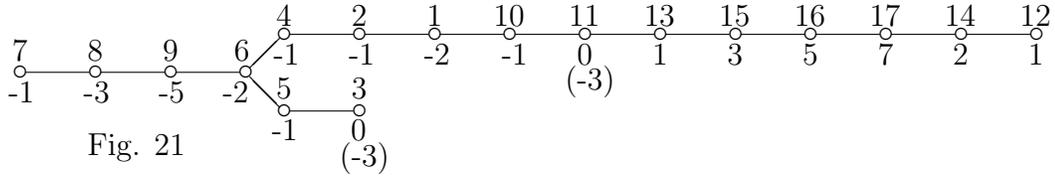
\setlength{\unitlength}{0.5cm}
\begin{picture}(10,3.5)
\put(1,-1.5){Fig. 21}
\multiput(-1,0.5)(2,0){4}{$\circ$}
\multiput(-0.65,0.7)(2,0){3}{\line(1,0){1.7}}
\put(-1.1,-0.1){-$1$}
\put(0.9,-0.1){-$3$}
\put(2.9,-0.1){-$5$}
\put(4.6,-0.1){-$2$}
\put(-1,1){$7$}
\put(1,1){$8$}
\put(3,1){$9$}
\put(4.9,1){$6$}

\put(5.32,0.57){\line(1,-1){0.79}}
\put(6.02,-0.53){$\circ$} 
\put(6.37,-0.33){\line(1,0){1.7}}
\put(8.02,-0.53){$\circ$}
\put(5.9,-1.1){-$1$}
\put(8.0,-1.1){$0$}
\put(7.7,-1.8){(-$3$)}
\put(6.02,0){$5$}
\put(8.02,0){$3$}

\put(5.32,0.82){\line(1,1){0.79}}
\multiput(6.02,1.51)(2,0){11}{$\circ$} 
\multiput(6.37,1.71)(2,0){10}{\line(1,0){1.7}}
\put(5.95,0.9){-$1$}
\put(7.95,0.9){-$1$}
\put(9.95,0.9){-$2$}
\put(12.02,0.9){-$1$}
\put(14.02,0.9){$0$}
\put(16.02,0.9){$1$}
\put(18.02,0.9){$3$}
\put(20.02,0.9){$5$}
\put(22.02,0.9){$7$}
\put(24.02,0.9){$2$}
\put(26.02,0.9){$1$}
\put(13.7,0.25){(-$3$)}

\put(6.02,1.95){$4$}
\put(8.02,1.95){$2$}
\put(10.02,1.95){$1$}
\put(11.8,1.95){$10$}
\put(13.8,1.95){$11$}
\put(15.8,1.95){$13$}
\put(17.8,1.95){$15$}
\put(19.8,1.95){$16$}
\put(21.8,1.95){$17$}
\put(23.8,1.95){$14$}
\put(25.8,1.95){$12$}

\end{picture}
\vskip 1.8cm

We will now create some branches from the (-5)-curve (creation number 9) and (-2)-curve (creation number 6) on the above graph. Specifically, from the (-5)-curve we will create 14 length 1 branches, 9 length 3 branches and 5 long forked branches, the same way as in the First Framework. They will be mapped to the branches on $Y$ in the same fashion as in the First Framework (so the total degree of the map is $28=14\cdot 2=9\cdot 3 +1\cdot 1= 1\cdot 23+ 5\cdot 1$). From the forked (-2)-curve we will create 16 branches of length 2, to be mapped 1-to-1 to the branch on $Y$ that ends with the 0-curve with self-intersection (-1). We will also construct 2 branches of length 4 and 4 branches of length 6, to be mapped to the branch on $Y$ that ends with the 0-curve with self-intersection (-2) with degrees 5 and 3 respectively (see below). Finally, we will modify the long branch from the (-2)-curve by contracting, in order, the curves numbered 1, 10, and 2, and blowing up the intersection of the (-2)-curve with the branch three times to get (-3)-curve, (-5)-curve, and (-7)-curve. This branch will be sent to  the branch on $Y$ that ends with the 0-curve with self-intersection (-1) with degree 7. Finally, like in the First Framework, we will break the edge between the two multi-forked curves, this time by creating 31 new curves, with the following $\bar{K}$ labels, from left to right: -92, -87, -82, -77, -72, -67, -62, -57, -52, -47, -42, -37, -69, -32, -27, -22, -17, -46, -29, -12, -7, -23, -16, -9, -11, -13, -15, -17, -19, -21, -23.

The following picture shows the general outline of the map, similar to Figure 10 for the First Framework.

\setlength{\unitlength}{0.5cm}
\begin{picture}(10,6)
\thicklines
\put(4.2,4.8){\bf Second Framework}
\multiput(-1,0.5)(2,0){2}{\color{Plum}$\circ$}
\multiput(3,0.5)(2,0){2}{\color{red}$\circ$}
\multiput(7,0.5)(2,0){2}{\color{green}$\circ$}

\multiput(-0.65,0.7)(2,0){2}{\color{Plum}\line(1,0){1.7}}
\multiput(5.35,0.7)(2,0){2}{\color{green}\line(1,0){1.7}}
\multiput(3.3,0.68)(0.2,0){9}{\color{Cyan}.}
\put(-0.9,-0.05){-$1$}
\put(0.92,-0.05){-$3$}
\put(2.75,-0.05){-$5$}
\put(4.43,-0.05){-$2$}
\put(6.92,-0.05){-$1$}
\put(8.95,-0.05){$0$}
\put(8.4,-0.65){(-$1$)}
\put(10,0.51){\circle{1.05}}
\put(9.47,0.34){\tiny $\times\! 16$}

{\color{blue}
\put(5.33,0.57){\line(1,-1){0.78}}
\put(6.02,-0.53){$\circ$} 
\put(6.37,-0.33){\line(1,0){1.7}}
\put(8.02,-0.53){$\circ$}
}
\put(5.95,-1.1){-$1$}
\put(7.9,-1.1){$0$}
\put(7.4,-1.7){(-$3$)}

{\color{blue}
\put(5.27,0.55){\line(1,-3){0.88}}
\multiput(6.02,-2.45)(2,0){4}{$\circ$} 
\multiput(6.37,-2.25)(2,0){3}{\line(1,0){1.7}}
}
\put(5.95,-3.1){-$5$}
\put(7.9,-3.1){-$3$}
\put(9.95,-3.1){-$1$}
\put(12.3,-2.9){$0$}
\put(12,-3.5){(-$1$)}
\put(12.4,-1.55){\circle{1}}
\put(12.,-1.69){\tiny $\times 2$}

{\color{blue}
\put(5.2,0.55){\line(1,-5){0.9}}
\multiput(5.93,-4.3)(2,0){6}{$\circ$} 
\multiput(6.28,-4.1)(2,0){5}{\line(1,0){1.7}}
}
\put(5.8,-4.9){-$3$}
\put(7.8,-4.9){-$1$}
\put(10,-4.9){$0$}
\put(12,-4.9){$1$}
\put(14.15,-4.9){$2$}
\put(16,-4.9){$3$}
\put(16.4,-5.3){$\nwarrow$}
\put(16.3,-5.8){\tiny $type\  3$}
\put(9.6,-5.6){(-$2$)}
\put(16.34,-3.35){\circle{1.2}}
\put(15.89,-3.54){\tiny $\times 4$}

{\color{green}
\put(5.33,0.85){\line(1,1){0.78}}
\multiput(6.02,1.51)(2,0){6}{$\circ$} 
\multiput(6.37,1.71)(2,0){5}{\line(1,0){1.7}}
\put(16.3,1.65){\line(1,-1){0.84}}
\multiput(17.02,0.5)(2,0){3}{$\circ$} 
\multiput(17.37,0.7)(2,0){2}{\line(1,0){1.7}}
\put(21.33,0.85){\line(1,1){0.78}}
\multiput(22.02,1.51)(2,0){2}{$\circ$} 
\multiput(22.37,1.71)(2,0){1}{\line(1,0){1.7}}
}
\put(5.9,1){-$7$}
\put(7.9,1){-$5$}
\put(10,1){-$3$}
\put(11.9,1){-$1$}
\put(14,1){$0$}
\put(16.,1){$1$}
\put(17,-0.1){$3$}
\put(19,-0.1){$5$}
\put(21,-0.1){$7$}
\put(22.1,1){$2$}
\put(24,1){$1$}
\put(13.6,0.4){(-$1$)}
\put(21.4,-0.3){$\nwarrow$}
\put(21.3,-0.8){\tiny $type\  3$}

{\color{Cyan}
\put(3.27,0.89){\line(1,3){0.88}}
\multiput(4.02,3.43)(2,0){4}{$\circ$}
\multiput(4.37,3.63)(2,0){4}{\line(1,0){1.7}}
}
\put(12.02,3.43){\color{red}$\circ$}
{\color{blue}
\multiput(14.02,3.43)(2,0){2}{$\circ$}
\multiput(12.37,3.63)(2,0){2}{\line(1,0){1.7}}
}
\put(4,2.9){-$4$}
\put(5.9,2.9){-$3$}
\put(7.9,2.9){-$2$}
\put(9.9,2.9){-$1$}
\put(11.7,2.9){-$2$}
\put(13.9,2.9){-$1$}
\put(16,2.9){$0$}
\put(15.7,2.3){(-$2$)}

{\color{green}
\put(12.32,3.8){\line(1,2){0.83}}
\put(13.02,5.4){$\circ$} 
\put(13.37,5.6){\line(1,0){1.7}}
\put(15.02,5.4){$\circ$}
}
\put(13.02,4.9){-$1$}
\put(15.02,4.9){$0$}
\put(14.78,4.2){(-$1$)}
\put(14.05,4.51){\circle{1}}
\put(13.6,4.32){\tiny $\times 5$}

{\color{OliveGreen}
\put(3.29,0.55){\line(1,-2){0.86}}
\put(4.01,-1.54){$\circ$} 
}
\put(3.6,-2.2){-$4$}
\put(3.4,-1.2){\circle{1.2}}
\put(2.85,-1.32){\tiny $\times \! 14$}

{\color{Plum}
\multiput(3.11,0.83)(-1,1){3}{\line(-1,1){0.75}}
\multiput(2.04,1.5)(-1,1){3}{$\circ$}
}
\put(1.8,.95){-$9$}
\put(0.8,1.95){-$4$}
\put(-0.45,2.95){-$3$}
\put(0.87,3.91){\circle{1}}
\put(0.44,3.72){\tiny $\times 9$}


\multiput(0.5,-9.52)(2,0){2}{\color{Plum} $\circ$}
\put(4.5,-9.52){\color{red} $\circ$}
\multiput(6.5,-9.52)(2,0){5}{\color{Cyan} $\circ$}
\put(16.5,-9.52){\color{red} $\circ$}
\multiput(18.5,-9.52)(2,0){2}{\color{blue}$\circ$}

\multiput(0.85,-9.32)(2,0){2}{\color{Plum} \line(1,0){1.7}}
\multiput(4.86,-9.32)(2,0){6}{\color{Cyan} \line(1,0){1.69}}
\multiput(16.85,-9.32)(2,0){2}{\color{blue}\line(1,0){1.7}}

\put(0.4,-10.1){-$1$}
\put(2.4,-10.1){-$3$}
\put(4.15,-10.1){-$5$}
\put(6.4,-10.1){-$4$}
\put(8.4,-10.1){-$3$}
\put(10.4,-10.1){-$2$}
\put(12.4,-10.1){-$1$}
\put(14.4,-10.1){-$1$}
\put(16.4,-10.1){-$2$}
\put(18.4,-10.1){-$1$}
\put(20.5,-10.1){$0$}
\put(20.2,-10.7){(-$2$)}

\put(4.75,-9.45){\color{OliveGreen} \line(1,-2){0.86}}
\put(5.47,-11.53){\color{OliveGreen}$\circ$} 
\put(5.4,-12.05){-$2$}

{\color{green}
\put(16.81,-9.14){\line(1,2){0.82}}
\put(17.5,-7.55){$\circ$} 
\put(17.86,-7.35){\line(1,0){1.69}}
\put(19.5,-7.55){$\circ$}
}
\put(17.6,-8.1){-$1$}
\put(19.5,-8.1){$0$}
\put(19.2,-8.7){(-$1$)}


{\color{Plum}
\qbezier(-0.9,0.57)(-1.5,-4)(0.46,-8.88)
\put(0.28,-9.2){$\triangleleft$}
\qbezier(0.25,3.55)(0.7,-6)(0.8,-8.68)
\put(0.5,-9.16){$\triangledown$}
}

{\color{red}
\qbezier(2.83,0.57)(0.6,-1)(4.2,-8.86)
\put(4.01,-9.18){$\triangleleft$}
}

{\color{red}
\qbezier(4.62,0.55)(4,-7)(9,-7)
\qbezier(9,-7)(12,-7)(15.63,-9.05)
\put(15.5,-9.25){\tiny $\triangle$}
}

{\color{OliveGreen}
\qbezier(3.42,-1.48)(4.5,-7)(4.92,-10.62)
\put(4.62,-11.14){$\triangledown$}
}

{\color{red}
\qbezier(11.25,3.53)(15.4,-7)(15.58,-8.63)
\put(15.3,-9.1){$\triangledown$}
}

{\color{blue}
\put(18.55,-9.26){\tiny $\triangle$}
\put(19.14,-9.1){\tiny $\triangledown$}
\put(19.5,-9.25){\tiny $\triangle$}
\qbezier(6.9,-0.46)(9,-5.5)(18.68,-9.06)
\qbezier(10.97,-2.2)(19,0)(19.3,-8.74)
\qbezier(8.88,-4.22)(13.2,-7.5)(18.68,-9.06)
\qbezier(14.98,3.54)(22.3,-5.5)(19.9,-9.05)

}

{\color{green}
\qbezier(7.67,0.55)(10,-5.4)(17.5,-7.05)
\put(17.4,-7.25){\tiny $\triangle$}
\qbezier(12.71,1.59)(17,-2)(18.04,-6.79)
\put(17.86,-7.12){\tiny $\triangledown$}
\qbezier(13.67,5.52)(23.5,-3)(18.59,-7.05)
\put(18.2,-7.25){\tiny $\triangle$}
}

{\color{Cyan}
\qbezier(2.3,0.5)(3,-6.5)(7,-8.69)
\put(6.8,-9){$\triangleleft$}
}

\put(-2.7,-4){\color{Plum}$deg\ 3$}
\put(-0.3,-5.5){\color{red}$deg\ 28$}
\put(1.9,-7){\color{OliveGreen}$deg\ 2$}
\put(3.91,-7.5){\color{Cyan}$deg\ 23$}
\put(7,-7.8){\color{red}$deg\ 23$}
\put(9.5,-6.8){\color{blue}$deg\ 3$}
\put(16.8,-5.5){\color{blue}$deg\ 5$}
\put(13.8,-1.5){\color{green}$deg\ 7$}

\end{picture}
\vskip 7cm

\setlength{\unitlength}{0.38cm}
\begin{picture}(10,1.5)
\thicklines
\put(6,2.5){Close-up of the (-5)...(-2) map (light blue arrow)}

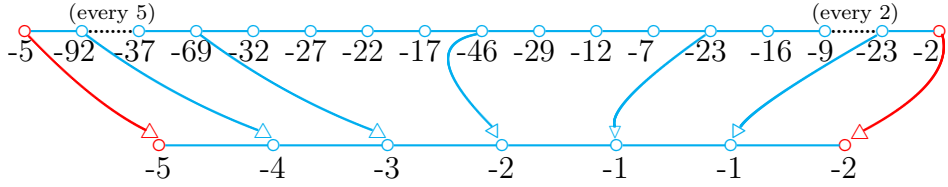
\multiput(-0.7,0.5)(32,0){12}{$\color{red} \circ$}
{\color{Cyan}
\multiput(1.3,0.5)(2,0){15}{$\circ$}
\multiput(-0.21,0.76)(2,0){1}{\line(1,0){1.6}}
\multiput(3.79,0.76)(2,0){12}{\line(1,0){1.6}}
\multiput(29.79,0.76)(2,0){1}{\line(1,0){1.6}}
}
\multiput(1.72,0.7)(0.2,0){8}{.}
\multiput(27.72,0.7)(0.2,0){8}{.}
\put(1.1,1.2){\tiny (every 5)}
\put(27.1,1.2){\tiny (every 2)}

\put(-1,-0.3){-$5$}
\put(0.6,-0.3){-$92$}
\put(2.7,-0.3){-$37$}
\put(4.7,-0.3){-$69$}
\put(6.7,-0.3){-$32$}
\put(8.7,-0.3){-$27$}
\put(10.7,-0.3){-$22$}
\put(12.7,-0.3){-$17$}
\put(14.7,-0.3){-$46$}
\put(16.7,-0.3){-$29$}
\put(18.7,-0.3){-$12$}
\put(20.7,-0.3){-$7$}
\put(22.7,-0.3){-$23$}
\put(25,-0.3){-$16$}
\put(27,-0.3){-$9$}
\put(28.7,-0.3){-$23$}
\put(30.7,-0.3){-$2$}

\multiput(4,-3.5)(24,0){2}{\color{red}$\circ$}
{\color{Cyan}
\multiput(8,-3.5)(4,0){5}{$\circ$}
\multiput(4.49,-3.24)(4,0){6}{\line(1,0){3.6}}
}
\put(3.8,-4.4){-$5$}
\put(7.8,-4.4){-$4$}
\put(11.8,-4.4){-$3$}
\put(15.8,-4.4){-$2$}
\put(19.8,-4.4){-$1$}
\put(23.8,-4.4){-$1$}
\put(27.8,-4.4){-$2$}

{\color{red}
\qbezier(-0.27,0.57)(2,-1.8)(3.8,-2.7)
\put(3.6,-3){\tiny $\triangle$}
}
{\color{Cyan}
\qbezier(1.38,0.61)(4.2,-1.5)(7.35,-2.7)
\put(7.2,-3){\tiny$\triangle$}
\qbezier(5.38,0.61)(8.2,-1.5)(11.35,-2.7)
\put(11.2,-3){\tiny$\triangle$}
\qbezier(15.01,0.65)(12.8,0.2)(15.54,-2.56)
\put(15.3,-3){\small $\triangleleft$}
\qbezier(23,0.61)(19.6,-1.8)(19.85,-2.58)
\put(19.6,-3){\tiny $\triangledown$}
\qbezier(29,0.6)(24.8,-1.8)(24.26,-2.58)
\put(23.95,-2.98){\small $\triangleright$}
}
{\color{red}
\qbezier(30.95,0.63)(31.4,-1)(28.22,-2.71)
\put(27.72,-3.02){\tiny$\triangle$}
}

\put(-1.5,4){Fig. 22}
\end{picture}
\vskip 2cm

We now look at the maps on the branches in more detail. We skip the branches from the (-5)-curve, as they are identical to those from the First Framework.
\begin{itemize}
\item For the length 4 branches from the (-2)-curve:

\setlength{\unitlength}{0.65cm}
\begin{picture}(10,3)
\put(-2.5,-1){Fig. 23}
\multiput(-0.98,1.51)(4,0){5}{$\circ$} 
\multiput(-0.73,1.65)(4,0){4}{\line(1,0){3.78}}
\put(-1.1,1.05){-$2$}
\put(2.9,1.05){-$5$}
\put(6.9,1.05){-$3$}
\put(10.9,1.05){-$1$}
\put(15.02,1.05){$0$}

\put(-1.5,2){(-$46$)}
\put(2.6,2){(-$1$)}
\put(6.6,2){(-$2$)}
\put(10.6,2){(-$3$)}
\put(14.6,2){(-$1$)}

\multiput(0,-2.14)(4,0){3}{$\circ$} 
\multiput(0.25,-2)(4,0){2}{\line(1,0){3.78}}
\put(-0.1,-2.6){-$2$}
\put(3.9,-2.6){-$1$}
\put(8,-2.6){$0$}
\put(-0.3,-3.2){(-$2$)}
\put(3.7,-3.2){(-$2$)}
\put(7.7,-3.2){(-$3$)}

\put(-0.07,-1.8){\small $\triangledown$}
\put(3.93,-1.8){\small $\triangledown$}
\put(5.36,-1.8){\small $\triangledown$}
\put(6.93,-1.8){\small $\triangledown$}
\put(8.2,-1.8){$\triangleright$}

\qbezier(-0.7,1.6)(-0.1,1)(0.15,-1.44)
\qbezier(3.3,1.6)(3.9,1)(4.15,-1.44)
\qbezier(7.1,1.6)(5.8,0.8)(5.57,-1.44)
\qbezier(11.1,1.6)(7.7,1)(7.14,-1.44)
\qbezier(15.1,1.6)(9.7,0)(8.36,-1.54)

\end{picture}
\vskip 2.5cm

\hskip -1cm $\phi_*(E_{-5})=F_{-1},$ $\phi_*(E_{-3})=\phi_*(E_{-1})=0,$  $\phi_*(E_{0})=F_{0},$ 

\hskip -1cm $\phi^*(F_{-1})=5E_{-5}+3E_{-3}+E_{-1},$ $\phi^*(F_{0})=5E_{0}+2E_{-1}+E_{-3}$ 

\item For the length 6 branches from the (-2)-curve:

\setlength{\unitlength}{0.65cm}
\begin{picture}(10,3.5)
\put(-2.2,-1){Fig. 24}
\multiput(-0.98,1.51)(2,0){7}{$\circ$} 
\multiput(-0.73,1.65)(2,0){6}{\line(1,0){1.78}}
\put(-1.1,1.05){-$2$}
\put(0.9,1.05){-$3$}
\put(2.9,1.05){-$1$}
\put(5.02,1.05){$0$}
\put(7.02,1.05){$1$}
\put(9.02,1.05){$2$}
\put(11.02,1.05){$3$}

\put(-1.5,2){(-$46$)}
\put(0.6,2){(-$1$)}
\put(2.6,2){(-$3$)}
\put(4.6,2){(-$2$)}
\put(6.6,2){(-$2$)}
\put(8.6,2){(-$2$)}
\put(10.6,2){(-$1$)}

\multiput(0,-2.14)(4,0){3}{$\circ$} 
\multiput(0.25,-2)(4,0){3}{\line(1,0){3.78}}
\put(-0.1,-2.6){-$2$}
\put(3.9,-2.6){-$1$}
\put(8,-2.6){$0$}
\put(3.7,-3.2){(-$2$)}
\put(7.7,-3.2){(-$3$)}

\put(-0.07,-1.8){\small $\triangledown$}
\put(3.76,-1.81){$\triangleleft$}
\put(6.06,-1.81){$\triangleleft$}
\put(7.93,-1.8){\small $\triangledown$}
\put(9.38,-1.8){\small $\triangledown$}
\put(10.8,-1.8){\small $\triangledown$}
\put(12,-1.95){\small $\triangledown$}

\qbezier(-0.7,1.6)(-0.1,1)(0.15,-1.44)
\qbezier(1.3,1.6)(2.5,0.8)(3.94,-1.56)
\qbezier(3.3,1.6)(4.7,0.8)(6.24,-1.56)
\qbezier(5.3,1.6)(8,0.8)(8.15,-1.45)
\qbezier(7.3,1.6)(9.58,1)(9.58,-1.45)
\qbezier(9.3,1.6)(11,1)(11,-1.45)
\qbezier(11.3,1.6)(12.2,1)(12.2,-1.6)

\end{picture}
\vskip 2.3cm
$\phi_*(E_{-3})=F_{-1},$   $\phi_*(E_{-1})=0,$ $\phi_*(E_{0})=F_{0},$

$\phi_*(E_{1})=\phi_*(E_{2})=0,$  $\phi (E_3)$ is a curve, generically in $\mathbb A^2,$

$\phi^*(F_{-1})=3E_{-3}+ E_{-1},$ $\phi^*(F_{0})=3E_{0}+ E_{-1}+(2E_1+E_2)$ 

\item For the long branch from the (-2)-curve:

\setlength{\unitlength}{.5cm}
\begin{picture}(10,2.8)
\put(-3,-2){Fig. 25}
\multiput(-2.5,0.6)(2,0){12}{$\circ$} 
\multiput(-2.15,0.8)(2,0){11}{\line(1,0){1.7}}
\put(-2.6,0){-$2$}
\put(-0.6,0){-$7$}
\put(1.4,0){-$5$}
\put(3.4,0){-$3$}
\put(5.4,0){-$1$}
\put(7.5,0){$0$}
\put(9.5,0){$1$}
\put(11.5,0){$3$}
\put(13.5,0){$5$}
\put(15.5,0){$7$}
\put(17.5,0){$2$}
\put(19.5,0){$1$}

\put(-3.1,1.1){(-$46$)}
\put(-1,1.1){(-$1$)}
\put(1,1.1){(-$2$)}
\put(3,1.1){(-$2$)}
\put(5,1.1){(-$3$)}
\put(7,1.1){(-$1$)}
\put(9,1.1){(-$3$)}
\put(11,1.1){(-$2$)}
\put(13,1.1){(-$2$)}
\put(15,1.1){(-$1$)}
\put(17,1.1){(-$4$)}
\put(19,1.1){(-$3$)}

\multiput(1,-3)(4,0){3}{$\circ$} 
\multiput(1.35,-2.8)(4,0){3}{\line(1,0){3.7}}
\put(0.9,-3.6){-$2$}
\put(4.9,-3.6){-$1$}
\put(9,-3.6){$0$}
\put(4.7,-4.3){(-$2$)}
\put(8.7,-4.3){(-$1$)}

\put(9,-0.8){$E_1$}
\put(19,-0.8){$E_1'$}

\put(0.8,-2.55){\small $\triangleleft$}
\put(4.8,-2.55){\small $\triangleleft$}
\put(6,-2.7){\tiny $\triangledown$}
\put(7,-2.7){\tiny $\triangledown$}
\put(8,-2.7){\tiny $\triangledown$}
\put(9,-2.55){\tiny $\triangledown$}
\put(10,-2.7){\tiny $\triangledown$}
\put(11,-2.7){\tiny $\triangledown$}
\put(12,-2.7){\tiny $\triangledown$}
\put(13,-2.7){\tiny $\triangledown$}

\qbezier(-2.15,0.73)(-1,0.5)(1,-2.24)
\qbezier(-0.15,0.73)(4,-1)(5,-2.24)
\qbezier(1.85,0.73)(6,-1)(6.22,-2.38)
\qbezier(3.85,0.73)(7,-1)(7.22,-2.38)
\qbezier(5.85,0.73)(8,-1)(8.22,-2.38)
\qbezier(7.85,0.73)(9,-1)(9.22,-2.23)
\qbezier(9.82,0.7)(10.2,0)(10.22,-2.38)
\qbezier(11.59,0.7)(11.25,0)(11.22,-2.38)
\qbezier(13.58,0.71)(12.25,0)(12.22,-2.38)
\qbezier(15.56,0.73)(13.25,-0.7)(13.22,-2.38)

\end{picture}
\vskip 2.5cm

$\phi_*(E_{-7})=F_{-1},$   $\phi_*(E_{-5})=\phi_*(E_{-3})=\phi_*(E_{-1})=0,$

$\phi_*(E_{0})=F_{0},$  $\phi_*(E_{1})=\phi_*(E_{3})=\phi_*(E_{5})=0,$ $\phi (E_7)$ is a curve, generically in $\mathbb A^2,$ $\phi (E_2)=\phi(E_1')$ is a point in $\mathbb A^2,$

$\phi^*(F_{-1})=7E_{-7}+ (5E_{-5}+3E_{-3}+E_{-1}),$

$\phi^*(F_{0})=7E_{0}+(3E_{-1}+2E_{-3}+E_{-5})+(3E_1+2E_3+E_5)$ 

\item Finally, for the chain between the (-5)-curve and the (-2)-curve:

\setlength{\unitlength}{0.38cm}
\hskip -1cm \begin{picture}(10,4)
\put(-1,-3){Fig. 26}
\multiput(-0.7,0.5)(32,0){12}{$\circ$}

\multiput(1.3,0.5)(2,0){15}{$\circ$}
\multiput(-0.21,0.76)(2,0){1}{\line(1,0){1.6}}
\multiput(3.79,0.76)(2,0){12}{\line(1,0){1.6}}
\multiput(29.79,0.76)(2,0){1}{\line(1,0){1.6}}

\multiput(1.72,0.7)(0.2,0){8}{.}
\multiput(27.72,0.7)(0.2,0){8}{.}
\put(-0.5,2.2){\small $\bar{K}$ every 5, all (-$2$)}
\put(23.5,2.2){\small $\bar{K}$ every 2, all (-$2$)}
\put(2.2,1.1){$\downarrow$}
\put(28.2,1.1){$\downarrow$}

\put(-1,-0.3){-$5$}
\put(0.6,-0.3){-$92$}
\put(2.7,-0.3){-$37$}
\put(4.65,-0.3){-$69$}
\put(6.7,-0.3){-$32$}
\put(8.7,-0.3){-$27$}
\put(10.7,-0.3){-$22$}
\put(12.7,-0.3){-$17$}
\put(14.7,-0.3){-$46$}
\put(16.7,-0.3){-$29$}
\put(18.7,-0.3){-$12$}
\put(20.7,-0.3){-$7$}
\put(22.7,-0.3){-$23$}
\put(25,-0.3){-$16$}
\put(27,-0.3){-$9$}
\put(28.7,-0.3){-$23$}
\put(30.7,-0.3){-$2$}
\put(12.4,-1.3){$E_{-17}'$}
\put(22.4,-1.3){$E_{-23}'$}
\put(28.4,-1.3){$E_{-23}$}

\put(-1.5,1.1){\small (-$56$)}
\put(0.6,1.1){\small (-$1$)}
\put(2.6,1.1){\small (-$3$)}
\put(4.6,1.1){\small (-$1$)}
\put(6.6,1.1){\small (-$3$)}
\put(8.6,1.1){\small (-$2$)}
\put(10.6,1.1){\small (-$2$)}
\put(12.6,1.1){\small (-$4$)}
\put(14.6,1.1){\small (-$1$)}
\put(16.6,1.1){\small (-$2$)}
\put(18.6,1.1){\small (-$3$)}
\put(20.6,1.1){\small (-$5$)}
\put(22.6,1.1){\small (-$1$)}
\put(24.6,1.1){\small (-$2$)}
\put(26.6,1.1){\small (-$3$)}
\put(28.6,1.1){\small (-$1$)}
\put(30.3,1.1){\small (-$46$)}

\multiput(4,-3.5)(24,0){2}{$\circ$}
\multiput(8,-3.5)(4,0){5}{$\circ$}
\multiput(4.49,-3.24)(4,0){6}{\line(1,0){3.6}}

\put(3.8,-4.4){-$5$}
\put(7.8,-4.4){-$4$}
\put(11.8,-4.4){-$3$}
\put(15.8,-4.4){-$2$}
\put(19.8,-4.4){-$1$}
\put(23.8,-4.4){-$1$}
\put(27.8,-4.4){-$2$}

\put(7.6,-5.4){(-$2$)}
\put(11.6,-5.4){(-$2$)}
\put(15.6,-5.4){(-$2$)}
\put(19.6,-5.4){(-$3$)}
\put(23.6,-5.4){(-$3$)}

\put(18.3,-2.75){$F_{-1}'$}
\put(22.4,-2.75){$F_{-1}$}

\qbezier(-0.29,0.59)(2,-1.8)(3.8,-2.7)
\put(3.6,-3){\tiny $\triangle$}
\qbezier(1.75,0.61)(4.2,-1.5)(7.75,-2.7)
\put(7.5,-3){\tiny$\triangle$}
\qbezier(5.75,0.61)(8.2,-1.5)(11.65,-2.7)
\put(11.5,-3){\tiny$\triangle$}
\qbezier(15.4,0.65)(13.5,0.2)(15.84,-2.56)
\put(15.6,-3){\small $\triangleleft$}
\qbezier(23.4,0.61)(20.2,-1)(20.2,-2.58)
\put(19.97,-3){\tiny $\triangledown$}
\qbezier(29.4,0.61)(25,-2)(24.56,-2.58)
\put(24.25,-2.98){\small $\triangleright$}

\qbezier(31.75,0.63)(32.4,-1)(28.82,-2.71)
\put(28.32,-3.02){\tiny$\triangle$}

\end{picture}
\vskip 2.3cm

On the above picture it is shown where the curves of type 1 go. The curves of type 2 go to the intersections of the ``neighboring" curves of type 1. Note that we have several pairs of curves with the same $\bar{K}$ label, including the curve $E_{-17}$, hidden between $E_{-9}$ and $E_{-23}$. For all curves of type 1 (not including the (-5)-curve and the (-2)-curve) we have $f=1$ and $e=23$. The $\phi^*$ is given by the following formulas. 

\hskip -1cm $\phi^* F_{-4}=23E_{-92}+(21E_{-87}+19E_{-82}+...+5E_{-47}+3E_{-42}+E_{-37}),$

\hskip -1cm $\phi^* F_{-3}=23E_{-69}+(11E_{-37}+10E_{-42}+...+3E_{-77}+2E_{-82}+E_{-87})+$

\hskip -1cm $(10E_{-32}+7E_{-27}+4E_{-22}+E_{-17}'),$

\hskip -1cm $\phi^* F_{-2}\!=\!\!23E_{-46}+\!(\!7E_{-17}'+\!5E_{-22}+\!3E_{-27}+\!E_{-32}\!)+\!(\!14E_{-29}+\!5E_{-2}+\!E_{-7}\!),$

\hskip -1cm $\phi^* F_{-1}'=23E_{-23}'+(5E_{-7}+2E_{-12}+E_{-29})+(15E_{-16}+7E_{-9}+6E_{-11}+$

\hskip -1cm $5E_{-13}+4E_{-15}+3E_{-17}+2E_{-19}+E_{-21}),$

\hskip -1cm $\phi^* F_{-1}=23E_{-23}+(20E_{-21}+17E_{-19}+14E_{-17})+(11E_{-15}+8E_{-3}+$

\hskip -1cm $5E_{-11}+2E_{-9}+E_{-16})$

\end{itemize}

Again, it is tedious but not hard to check that the projection formula is true. Like in the First Framework, we can also calculate the degrees of the corresponding polynomials. With the same choice of coordinates on $Y$ and the initial blowup, we get that the valuations of $y_1$ and $y_2$ on the  curves numbered 12, 13, and 14 on Figure 20 are (-15) and (-10) respectively.  From the formulas below Figure 25, $\phi^*(y_1)$ and $\phi^*(y_2)$  have poles of orders 60 and 40 respectively on $E_{-1}$ and 105 and 70 on $E_0$. These are the curves numbered 4 and 11 on Figure 21. Reconstructing  the curves numbered 2, 10, and 1 on Figure 21, we get the following orders of poles: (165,110), (270,180), and, finally, $(435,290)$. So with the suitable choice of coordinates on the source and the target planes, our map $\phi$ should be given by a pair of polynomials of degrees $(435,290)$. We conclude our discussion of the Second Framework by the dessins for the two red arrow Belyi maps.
\vskip 0.5cm

{\bf Belyi map of (-5)-curves, the Second Framework}

The degree of the map is 28. The map is ramified above three points: intersection with the (-2)-curve, (-3)-curve, and (-4)-curve, that we will identify with $\{0\}$, $\{\infty\},$ and $\{1\}$ respectively. Above $\{0\}$, we have 14 ramification points of order 2, so the corresponding dessin is ``clean". Above $\{\infty\}$ we have 9 points of order 3 and 1 point of order 1, and above $\{1\}$ we have one point of order 23 and 5 points of order 1. A simple combinatorial analysis leads to the following clean dessin d'enfant (not unique, there are some options):

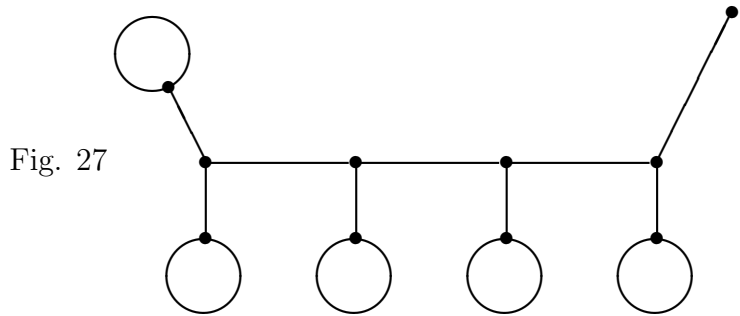
\setlength{\unitlength}{1cm}
\begin{picture}(10,2.5)
\put(0.5,0){Fig. 27}
\thicklines
\multiput(3,0)(2,0){4}{$\bullet$}
\multiput(3.15,0.1)(2,0){3}{\line(1,0){1.86}}

\multiput(3,-1)(2,0){4}{$\bullet$}
\put(10,2){$\bullet$}
\put(9.13,0.17){\line(1,2){0.93}}

\put(2.51,1){$\bullet$}
\put(3.1,0.15){\line(-1,2){0.45}}
\put(2.4,1.55){\circle{1}}

\multiput(3.08,-1.4)(2,0){4}{\circle{1}}
\multiput(3.11,-0.9)(2,0){4}{\line(0,1){1}}
\end{picture}
\vskip 2.5cm

{\bf Belyi map of (-2)-curves, the Second Framework}

The degree of the map is 23. The map is ramified above three points. They correspond to the branches that end with the 0-curve with the self-intersection  (-2), with the 0-curve with the self-intersection (-1), and with the forked (-5)-curve. We will identify them as $\{0\}$,  $\{1\}$, and $\{\infty\}$ respectively. On the (-2)-curve on $Z$ we will identify with $\{\infty\}$ the unique point that is sent to $\{\infty\}$, so that the Belyi map is given by a polynomial. Then above $\{0\}$ we have one point of order 1,  4 points of order 3, and 2 points of order 5. Above $\{1\}$ we have one point of order 7 and 16 points of order 1. A corresponding dessin d'enfant (unique as a graph, but  not as a dessin) is the following:

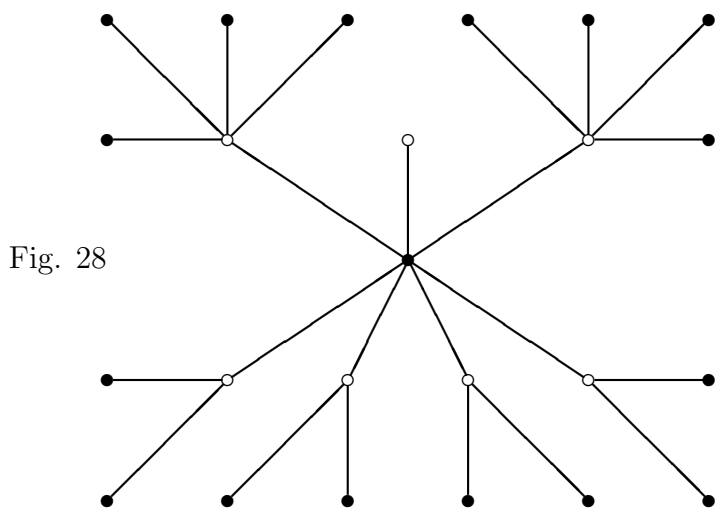
\setlength{\unitlength}{0.8 cm}
\begin{picture}(10,5)
\put(0.5,0){Fig. 28}
\thicklines
\multiput(2,4)(2,0){6}{$\bullet$}
\multiput(2,2)(10,0){2}{$\bullet$}
\multiput(4,2)(3,0){3}{$\circ$}
\multiput(4,-2)(2,0){4}{$\circ$}
\multiput(2,-2)(10,0){2}{$\bullet$}
\multiput(2,-4)(2,0){6}{$\bullet$}
\put(7,0){$\bullet$}

\multiput(7.05,0.19)(3,-2){2}{\line(-3,2){2.83}}
\multiput(7.21,0.19)(-3,-2){2}{\line(3,2){2.83}}
\put(7.13,0.19){\line(0,1){1.85}}
\put(7.12,0.1){\line(-1,-2){0.94}}
\put(7.16,0.1){\line(1,-2){0.94}}

\multiput(2.23,2.14)(8,0){2}{\line(1,0){1.82}}
\multiput(2.23,-1.86)(8,0){2}{\line(1,0){1.82}}
\multiput(2.21,4.06)(6,0){2}{\line(1,-1){1.85}}
\multiput(8.21,-1.94)(2,0){2}{\line(1,-1){1.85}}
\multiput(4.21,2.21)(6,0){2}{\line(1,1){1.85}}
\multiput(2.21,-3.79)(2,0){2}{\line(1,1){1.85}}
\multiput(4.13,2.23)(6,0){2}{\line(0,1){1.82}}
\multiput(6.13,-3.77)(2,0){2}{\line(0,1){1.82}}
\end{picture}
\vskip 4cm

\section{More Frameworks}

We will first construct several frameworks that are closely related to the First Framework. In fact, they will have the same target graph. We will call these frameworks ``isotopes" of the First Framework.

Besides the length 3 branches from (-2)-curve that are sent down with degree 3, as in Picture 12, there can also be branches that are sent down with degree 1, with a curve of type 3 at the end of the branch:

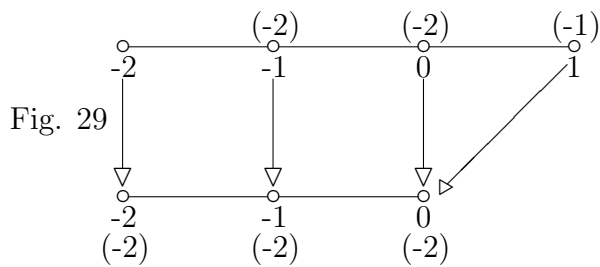
\setlength{\unitlength}{1cm}
\begin{picture}(10,1.5)
\put(-0.4,-0.5){Fig. 29}
\multiput(1,0.5)(2,0){4}{$\circ$}
\multiput(1.17,0.595)(2,0){3}{\line(1,0){1.86}}
\multiput(1,-1.5)(2,0){3}{$\circ$}
\multiput(1.17,-1.4)(2,0){2}{\line(1,0){1.86}}

\put(0.95,0.2){-$2$}
\put(2.95,0.2){-$1$}
\put(5,0.2){$0$}
\put(7,0.2){$1$}
\put(2.8,0.75){(-$2$)}
\put(4.8,0.75){(-$2$)}
\put(6.8,0.75){(-$1$)}

\put(0.95,-1.8){-$2$}
\put(2.95,-1.8){-$1$}
\put(5,-1.8){$0$}
\put(0.8,-2.2){(-$2$)}
\put(2.8,-2.2){(-$2$)}
\put(4.8,-2.2){(-$2$)}

\put(0.95,-1.26){$\triangledown$}
\put(2.95,-1.26){$\triangledown$}
\put(4.95,-1.26){$\triangledown$}
\put(5.3,-1.37){$\triangleright$}

\put(1.1,0.17){\line(0,-1){1.2}}
\put(3.1,0.17){\line(0,-1){1.2}}
\put(5.1,0.17){\line(0,-1){1.2}}
\put(7.01,0.36){\line(-1,-1){1.59}}

\end{picture}
\vskip 2.3cm

Here $\phi (E_1)$ intersects $F_0$ transversally at one point, all other maps are 1-to-1.

Additionally, the map on Figure 13 can be generalized to odd ramification $e=2k+1$ higher than 5, as follows.

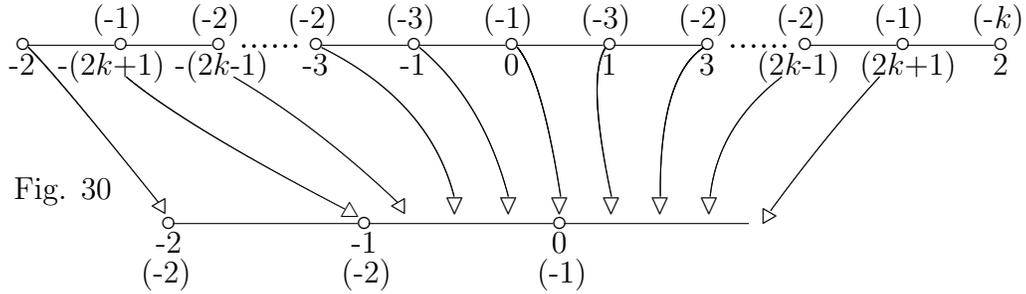
\setlength{\unitlength}{0.65cm}
\begin{picture}(10,3.5)
\put(-1,-1.5){Fig. 30}
\multiput(-0.98,1.51)(2,0){11}{$\circ$} 
\multiput(-0.73,1.65)(2,0){2}{\line(1,0){1.78}}
\multiput(5.27,1.65)(2,0){4}{\line(1,0){1.78}}
\multiput(15.27,1.65)(2,0){2}{\line(1,0){1.78}}
\put(3.6,1.6){\large ......}
\put(13.6,1.6){\large ......}

\put(-1.1,1.05){-$2$}
\put(-0.1,1.05){-(2$k$+1)}
\put(2.3,1.05){-(2$k$-1)}
\put(4.9,1.05){-$3$}
\put(6.9,1.05){-$1$}
\put(9.02,1.05){$0$}
\put(11.02,1.05){$1$}
\put(13.02,1.05){$3$}
\put(14.2,1.05){(2$k$-1)}
\put(16.3,1.05){(2$k$+1)}
\put(19.02,1.05){$2$}

\put(0.6,2){(-$1$)}
\put(2.6,2){(-$2$)}
\put(4.6,2){(-$2$)}
\put(6.6,2){(-$3$)}
\put(8.6,2){(-$1$)}
\put(10.6,2){(-$3$)}
\put(12.6,2){(-$2$)}
\put(14.6,2){(-$2$)}
\put(16.6,2){(-$1$)}
\put(18.6,2){(-$k$)}

\multiput(2,-2.14)(4,0){3}{$\circ$} 
\multiput(2.25,-2)(4,0){3}{\line(1,0){3.78}}
\put(1.9,-2.6){-$2$}
\put(5.9,-2.6){-$1$}
\put(10,-2.6){$0$}
\put(1.6,-3.2){(-$2$)}
\put(5.7,-3.2){(-$2$)}
\put(9.7,-3.2){(-$1$)}

\put(1.8,-1.8){$\triangleleft$}
\put(5.66,-1.87){\tiny $\triangle$}
\put(6.7,-1.81){$\triangleleft$}
\put(7.8,-1.8){\small $\triangledown$}
\put(8.9,-1.8){\small $\triangledown$}
\put(9.94,-1.8){\small $\triangledown$}
\put(11,-1.8){\small $\triangledown$}
\put(12,-1.8){\small $\triangledown$}
\put(13,-1.8){\small $\triangledown$}
\put(14.3,-2){$\triangleright$}

\qbezier(-0.7,1.6)(-0.1,1)(1.95,-1.56)
\qbezier(1.3,1)(2.5,0)(5.76,-1.7)
\qbezier(3.5,1)(5.2,0)(6.85,-1.57)
\qbezier(5.3,1.6)(7.5,0.5)(8,-1.44)
\qbezier(7.3,1.6)(8.6,0.5)(9.11,-1.44)
\qbezier(9.3,1.6)(9.8,1)(10.18,-1.44)
\qbezier(11.1,1.6)(10.7,1)(11.22,-1.44)
\qbezier(13.1,1.6)(12.2,1)(12.22,-1.44)
\qbezier(14.7,1)(13.4,0)(13.24,-1.44)
\qbezier(16.7,1)(15.6,-0.2)(14.48,-1.73)

\end{picture}
\vskip 3cm
$\phi_*(E_{-(2k+1)})=F_{-1},$   $\phi_*(E_{-(2k-1)})=...=\phi_*(E_{-1})=0,$ $\phi_*(E_{0})=F_{0},$

$\phi_*(E_{1})=...=\phi_*(E_{(2k-1)})=0,$  $\phi (E_{(2k+1)})$ is a curve, generically in $\mathbb A^2,$ intersecting $F_0$ at one point; $\phi(E_2)$ is a point in $\mathbb A^2,$

$\phi^*(F_{-1})=(2k+1)E_{-5}+(2k-1)E_{-(2k-1)}+...+ 3E_{-3}+E_{-1},$

$\phi^*(F_{0})=(2k+1)E_{0}+ (E_{-(2k-1)}+2E_{-(2k-3)}+...+(k-1)E_{-3}+kE_{-1})+$

$(kE_1+(k-1)E_3+...+2E_{(2k-3)}+E_{(2k-1)})$ 

This gives us some options for the Belyi map between the (-2)-curves. Specifically, for every $k\in \{2,3,4,5,6\}$ we can have a Belyi map with the following ramification data:
\begin{itemize}
\item above $\{\infty\}$: 1 point with ramification 13 (this point can be chosen to be $\{\infty \}$);

\item above $\{0\}$: $(6-k)$ points with ramification 3 and $(3k-5)$ points with ramification 1.

\item above $\{1\}$: 1 point with ramification $(2k+1)$ and $(12-2k)$ points with ramification 1.
\end{itemize}

Note that when $k=2$ we get the Belyi map described in Figure 16. When $k=6$ this map can be given by the polynomial function $t\mapsto t^{13}+1$. It is not hard to draw the possible dessins for all $k$, this is left to the reader as a pleasant exercise. The corresponding framework is similar to the  First Framework, but the long branch from the (-2)-curve on $Z$ is replaced by the long branch from Figure 30, the number of green branches of length 2 is $(12-2k)$ instead of 8, the number of blue branches of length 3 is $(6-k)$ instead of 4, and there are also $(3k-5)$ blue branches like in Figure 29.

It is not hard to calculate the degrees of the corresponding polynomials. Wih the same convention for $(x_1,x_2)$ and $(y_1,y_2)$ as for the First Framework, we get that
$y_1$ has degrees $9k+9$ and $27k+18$ in $x_1$ and $x_2$ respectively, while $y_2$ has degrees $6k+6$ and $18k+12$ respectively. The total degrees are $(36k+27,24k+18)$, that is the following:
$$(99,66),\ (135,90),\ (171,114),\ (207,138),\ (243,162) $$

Finally, here is a sketch of a more complicated framework. It also has some commonalities with the First Framework, but it has three rational Belyi maps instead of two. Interestingly, it has no curves of type 4.

\setlength{\unitlength}{0.65cm}
\begin{picture}(10,6)
\thicklines
\put(4.5,5){\bf Three-dessin Framework}
\multiput(-1,0.5)(2,0){2}{\color{Plum}$\circ$}
\multiput(3,0.5)(2,0){2}{\color{red}$\circ$}
\multiput(6,0.5)(1,0){2}{\color{green}$\circ$}
\multiput(-0.75,0.64)(2,0){2}{\color{Plum}\line(1,0){1.78}}
\multiput(5.27,0.64)(1,0){3}{\color{green}\line(1,0){0.78}}
\multiput(3.26,0.62)(0.2,0){9}{\color{Cyan}.}
\put(8,0.5){\color{red} $\circ$}
{\color{Dandelion}
\put(8.24,0.74){\line(1,2){0.6}}
\put(8.75,1.9){$\circ$}
}
\put(-0.9,0.05){-$1$}
\put(0.92,0.05){-$3$}
\put(2.75,0.05){-$5$}
\put(4.6,0.05){-$2$}
\put(5.92,0.05){-$1$}
\put(6.92,0.05){-$1$}
\put(7.82,0.05){-$1$}
\put(9.05,0.05){$0$}
\put(8.91,-0.48){(-$3$)}
\put(9.98,0.69){\circle{1}}
\put(9.56,0.48){$\times 8$}

{\color{blue}
\put(5.25,0.55){\line(1,-1){0.84}}
\put(6.02,-0.53){$\circ$} 
\put(6.27,-0.38){\line(1,0){1.78}}
\put(8.02,-0.53){$\circ$}
}
\put(5.95,-0.97){-$1$}
\put(8.02,-0.97){$0$}
\put(7.78,-1.5){(-$2$)}

{\color{blue}
\put(5.23,0.53){\line(1,-3){0.9}}
\multiput(6.02,-2.45)(2,0){3}{$\circ$} 
\multiput(6.27,-2.31)(2,0){2}{\line(1,0){1.78}}
}
\put(5.95,-2.95){-$3$}
\put(7.95,-2.95){-$1$}
\put(9.97,-2.95){$0$}
\put(9.6,-3.48){(-$1$)}
\put(9.1,-2.95){\circle{1}}
\put(8.65,-3.14){$\times 4$}

{\color{green}
\put(5.25,0.75){\line(1,2){0.9}}
\multiput(6.02,2.51)(2,0){4}{$\circ$}
\multiput(6.27,2.65)(2,0){4}{\line(1,0){1.78}}
}
\put(14.02,2.51){\color{red} $\circ$}
{\color{Dandelion}
\multiput(16.02,2.51)(2,0){3}{$\circ$} 
\multiput(14.27,2.65)(2,0){3}{\line(1,0){1.78}}
\put(14.23,2.75){\line(1,2){0.8}}
\put(14.93,4.3){$\circ$}
}
\put(15,3.8){$0$}
\put(14.7,3.3){(-$1$)}
\put(14.1,3.79){\circle{1}}
\put(13.65,3.6){$\times 3$}

\put(5.95,2.05){-$5$}
\put(7.95,2.05){-$3$}
\put(10,2.05){-$4$}
\put(12.02,2.05){-$5$}
\put(13.6,2.05){-$1$}
\put(16.02,2.05){$0$}
\put(18.02,2.05){$1$}
\put(20.02,2.05){$2$}
\put(15.78,1.52){(-$1$)}
\put(18.65,1){\tiny $type\  3$}
\put(19.25,1.6){$\nearrow$}
\put(19,0.2){$\downarrow$}

{\color{Cyan}
\put(3.23,0.78){\line(1,3){0.9}}
\multiput(4.02,3.43)(2,0){2}{$\circ$}
\multiput(4.27,3.57)(2,0){2}{\line(1,0){1.78}}
}
\put(8.26,3.54){{\large ...} 1-to-1 to Y}
\put(4.02,2.95){-$4$}
\put(5.95,2.95){-$3$}
\put(5.25,2.81){\circle{1}}
\put(4.8,2.6){$\times 3$}

{\color{OliveGreen}
\put(3.23,0.53){\line(1,-2){0.89}}
\put(4.02,-1.53){$\circ$} 
}
\put(3.95,-2){-$4$}
\put(3.5,-1.45){\circle{1}}
\put(3.05,-1.67){$\times 8$}

{\color{Plum}
\multiput(3.15,0.75)(-1,1){3}{\line(-1,1){0.84}}
\multiput(2.06,1.51)(-1,1){3}{$\circ$}
}
\put(1.9,1.05){-$9$}
\put(0.9,2.05){-$4$}
\put(-0.3,3.05){-$3$}
\put(0.87,3.91){\circle{1}}
\put(0.42,3.7){$\times 5$}

{\color{Sepia}
\put(14.23,2.53){\line(1,-3){0.9}}
\multiput(15.02,-0.45)(2,0){3}{$\circ$} 
\multiput(15.27,-0.31)(2,0){2}{\line(1,0){1.78}}
\put(14.29,2.56){\line(1,-1){1}}
\put(15.22,1.28){$\circ$} 
}
\put(14.95,-0.95){$0$}
\put(16.95,-0.95){$1$}
\put(18.97,-0.95){$2$}
\put(14.73,-1.48){(-$2$)}
\put(15.22,0.83){$0$}
\put(15,0.32){(-$1$)}

{\color{Sepia}
\put(8.27,0.64){\line(1,0){0.78}}
\put(9,0.5){$\circ$}
}
\put(8.78,1.47){$0$}
\put(8.5,0.97){(-$1$)}


\multiput(0.5,-9.52)(2,0){2}{\color{Plum} $\circ$}
\put(4.5,-9.52){\color{red} $\circ$}
\multiput(6.5,-9.52)(2,0){4}{\color{Cyan} $\circ$}
\put(14.5,-9.52){\color{red} $\circ$}
\multiput(16.5,-9.52)(2,0){2}{\color{blue}$\circ$}

\multiput(0.75,-9.37)(2,0){2}{\color{Plum} \line(1,0){1.78}}
\multiput(4.75,-9.37)(2,0){5}{\color{Cyan}\line(1,0){1.78}}
\multiput(14.75,-9.37)(2,0){2}{\color{blue}\line(1,0){1.78}}

\put(0.4,-10){-$1$}
\put(2.4,-10){-$3$}
\put(4.15,-10){-$5$}
\put(6.4,-10){-$4$}
\put(8.4,-10){-$3$}
\put(10.4,-10){-$2$}
\put(12.4,-10){-$1$}
\put(14.4,-10){-$2$}
\put(16.4,-10){-$1$}
\put(18.5,-10){$0$}
\put(18.2,-10.6){(-$2$)}

\put(4.7,-9.48){\color{OliveGreen} \line(1,-2){0.88}}
\put(5.47,-11.53){\color{OliveGreen}$\circ$} 
\put(5.4,-12){-$2$}

{\color{green}
\multiput(14.72,-9.28)(0.5,1){2}{\line(1,2){0.4}}
\put(15.5,-7.55){$\circ$} 
\put(15,-8.55){$\circ$}
\put(15.75,-7.4){\line(1,0){0.81}}
}
\put(16.75,-7.4){\color{Sepia} \line(1,0){0.81}}
\put(17.5,-7.55){\color{Sepia} $\circ$}

\put(16.5,-7.55){\color{red} $\circ$}
{\color{Dandelion}
\put(16.72,-7.28){\line(2,3){1.5}}
\put(18.12,-5.09){$\circ$}
}

\put(15,-9.03){-$1$}
\put(15.5,-8.03){-$1$}
\put(16.5,-8.03){-$1$}
\put(17.5,-8.03){$0$}
\put(17.25,-8.55){(-$3$)}
\put(18.2,-5.5){$0$}
\put(17.95,-6.02){(-$1$)}


{\color{Plum}
\qbezier(-0.9,0.53)(-1.5,-2)(0.42,-8.93)
\put(0.28,-9.2){$\triangleleft$}
\qbezier(0.25,3.55)(0.5,-6)(0.74,-8.78)
\put(0.5,-9.16){$\triangledown$}
}

{\color{red}
\qbezier(2.82,0.55)(0.6,-1)(4.25,-8.93)
\put(4.08,-9.2){$\triangleleft$}
}

{\color{red}
\qbezier(4.7,0.55)(5,-6)(13.78,-9.05)
\put(13.65,-9.25){\small $\triangle$}
}

{\color{OliveGreen}
\qbezier(3.64,-1.42)(5,-2)(5.04,-10.83)
\put(4.81,-11.2){$\triangledown$}
}

{\color{blue}
\qbezier(9.35,-2.38)(13,-8)(17.25,-9.05)
\put(17.15,-9.25){\small $\triangle$}
\qbezier(7.45,-0.45)(18.2,-4.5)(17.83,-8.8)
\put(17.6,-9.2){$\triangledown$}
}

{\color{green}
\qbezier(5.1,2.53)(6,-6)(13.65,-8.1)
\put(13.5,-8.3){\small $\triangle$}
\qbezier(11.02,2.56)(10,-3)(14.12,-7.06)
\put(14,-7.29){\small $\triangle$}
}

{\color{Cyan}
\qbezier(2.9,0)(5,-5)(8.15,-8.95)
\put(8,-9.2){$\triangleleft$}
}

{\color{red}
\qbezier(12.65,2.55)(12.2,-5)(14.72,-7.06)
\put(14.6,-7.29){\small $\triangle$}
\qbezier(6.75,0.55)(8.8,-3)(14.72,-7.06)
\put(14.6,-7.29){\small $\triangle$}
}

{\color{Sepia}
\qbezier(13.39,-0.36)(12,-2)(15.6,-7.11)
\put(15.5,-7.29){\tiny $\triangle$}
\qbezier(13.79,1.45)(15,1)(15.98,-6.9)
\put(15.77,-7.27){\small $\triangledown$}
}

{\color{Dandelion}
\qbezier(13.32,4.5)(17.2,4)(16.38,-4.4)
\qbezier(14.37,2.56)(16,1)(16.35,-4.4)
\put(16.15,-4.75){\small $\triangledown$}
}

\put(-2.47,-3){\color{Plum}$deg\ 3$}
\put(0.6,-7){\color{red}$deg\ 16$}
\put(2.7,-8){\color{OliveGreen}$deg\ 2$}
\put(5.6,-8){\color{Cyan}$deg\ 13$}
\put(5.9,-6.4){\color{red}$deg\ 13$}
\put(8.99,-4.9){\color{blue}$deg\ 3$}
\put(6.22,-5){\color{green}$deg\ 5$}
\put(9.55,-3.2){\color{green}$deg\ 5$}
\put(11.15,-1.5){\color{red}$deg\ 5$}
\put(13.16,-5.5){\color{Sepia}$deg\ 2$}
\put(14.02,-2.5){\color{Sepia}$deg\ 3$}
\put(14.45,0.5){\color{Dandelion}$deg\ 2$}

\end{picture}
\vskip 7.8cm

\setlength{\unitlength}{0.38cm}
\begin{picture}(10,2)
\thicklines
\put(10.5,2){Close-up of the (-5)...(-2) map (light blue arrow)}
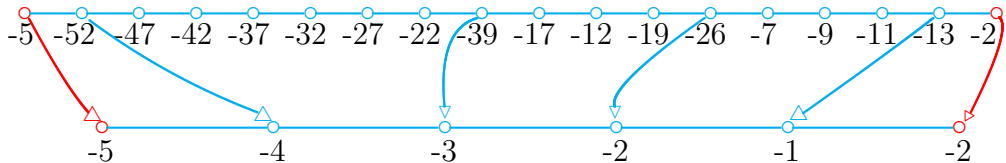
\multiput(-0.7,0.5)(34,0){12}{$\color{red} \circ$}
{\color{Cyan}
\multiput(1.3,0.5)(2,0){16}{$\circ$}
\multiput(-0.27,0.76)(2,0){17}{\line(1,0){1.62}}
}
\put(-1,-0.3){-$5$}
\put(0.6,-0.3){-$52$}
\put(2.7,-0.3){-$47$}
\put(4.7,-0.3){-$42$}
\put(6.7,-0.3){-$37$}
\put(8.7,-0.3){-$32$}
\put(10.7,-0.3){-$27$}
\put(12.7,-0.3){-$22$}
\put(14.7,-0.3){-$39$}
\put(16.7,-0.3){-$17$}
\put(18.7,-0.3){-$12$}
\put(20.7,-0.3){-$19$}
\put(22.7,-0.3){-$26$}
\put(25,-0.3){-$7$}
\put(27,-0.3){-$9$}
\put(28.7,-0.3){-$11$}
\put(30.7,-0.3){-$13$}
\put(32.7,-0.3){-$2$}

\multiput(2,-3.5)(30,0){2}{\color{red}$\circ$}
{\color{Cyan}
\multiput(8,-3.5)(6,0){4}{$\circ$}
\multiput(2.43,-3.24)(6,0){5}{\line(1,0){5.62}}
}
\put(1.8,-4.4){-$5$}
\put(7.8,-4.4){-$4$}
\put(13.8,-4.4){-$3$}
\put(19.8,-4.4){-$2$}
\put(25.8,-4.4){-$1$}
\put(31.8,-4.4){-$2$}

{\color{red}
\qbezier(-0.27,0.57)(1,-1.8)(1.8,-2.7)
\put(1.6,-3){\tiny $\triangle$}
}
{\color{Cyan}
\qbezier(1.48,0.52)(4.2,-1.5)(7.35,-2.7)
\put(7.2,-3){\tiny$\triangle$}
\qbezier(15.05,0.6)(13.65,0.2)(13.9,-2.5)
\put(13.65,-2.9){\tiny$\triangledown$}
\qbezier(23,0.6)(19.6,-1.8)(19.85,-2.5)
\put(19.6,-2.9){\tiny$\triangledown$}
\qbezier(31,0.6)(29.6,-0.5)(26.45,-2.73)
\put(25.95,-2.98){\tiny$\triangle$}
}
{\color{red}
\qbezier(32.95,0.63)(33.4,0.2)(31.95,-2.7)
\put(31.72,-3.02){\tiny$\triangleright$}
}

\put(0.5,3){Fig. 31}
\end{picture}
\vskip 2cm

It is a routine exercise to write down the exact maps between the branches. Most of them are the same as for the First Framework, and the rest are left to the reader.  This framework has three rational Belyi maps, above the three forked vertices on the graph of $Y$. Two of them, for the (-5)-curves and the (-2)-curves are the same as those in the First Framework.  The third one, for the (-1)-curves, has degree 5 and the following ramification data:

\begin{itemize}
\item above $\{\infty\}$: 1 point with ramification 5 (this point can be chosen to be $\{\infty \}$);

\item above $\{0\}$: 1 point with ramification 3 and 1 point with ramification 2;

\item above $\{1\}$: 1 point with ramification 2 and 3 points with ramification 1.
\end{itemize}

This Belyi map can be given by the function $t \mapsto \frac{1}{108}x^3 (x-5)^2$.  

Like in the previously considered frameworks, it is not hard to figure out the pair of degrees of the possible Keller map: $(108,72)$.

\section{Explanations and Comments}

As you have undoubtedly realized, the above frameworks are far from random. So I will first describe where they came from, and then briefly discuss some open questions that may help us solve the Jacobian Conjecture in dimension 2. Some of the proofs in this section are only sketched, and certain degree of familiarity with the  papers \cite{AmpleRamification} and \cite{DivisorialValuations} is required.

A lot of the discussion in this section will revolve around the curves with $\bar{K}$ label 0. They seem to be of great significance for the problem. We already know that for all other curves the $\bar{K}$ labels are sufficient to reconstruct their self-intersection numbers. It is also clear from the Keller Map Adjunction Formula that 0-curve on $Z$ is either of type 2 or  is sent by $\phi$ to a 0-curve on $Y$.

As noted in the Introduction, our frameworks basically come from repeatedly blowing up a point $\phi(\pi_*^{-1}(\infty))$ on $Y$ until $\pi_*^{-1}(\infty)$ becomes of type 1. I say ``basically", because in some cases it more natural to use not the line at infinity at $\mathbb P^2$ but to a Hirzebruch curve, which is a curve defined as follows:

\begin{Definition}  A Hirzebruch curve is the exceptional section of some Hirzebruch surface $F_n$, considered as a compactification of $\mathbb A^2$. This is a curve that is obtained after one blowup at infinity from $\mathbb P^2$. 
\end{Definition}

The Hirzebruch curve has $\bar{K}$ label (-1) and determinant label 0 (to be described below). In our first two frameworks (Figures 10 and 22), this is the curve with the $\bar{K}$ label (-1) on the long branch from the (-2)-curve. It is because we use this curve and not $\pi^{-1}_*(\infty)$ that at the end of the construction of $Z$ we contracted a couple of curves on $Z$. And we also stopped a little short of making this a type 1 curve, we just went until all curves of type 3 do not go through its image. In the Three-dessin Framework (Figure 31) this is the curve on $Z$ on which the third Belyi map is defined, but there we  did use $\pi^{-1}_*(\infty)$ and then contracted a couple of curves.  

Suppose $E$ on $Z$ is a Hirzebruch curve. Then it cannot be of type 3 or 4, because then a generic fiber on the Hirzebruch surface would be mapped entirely into $\mathbb A^2,$ which is clearly impossible. So the  procedure described above makes sense and must terminate.

Since we are resolving one divisorial valuation on $Y$, the graph of $Y$ will only have vertices (curves) with 1, 2, or 3 neighbors (i.e. with valency 1, 2, or 3). One can show, using the Keller condition and some adjunction inequalities, that all curves above the curves with valency 1 also have valency 1, and $f=1$ (cf. \cite{JC5}, Theorem 4.7). Likewise, all curves above curves with valency 2 must have valency 2. And above curves with valency 3 we have rational Belyi maps.

From the rules for the $\bar{K}$ labels, we can easily deduce that if in the process of constructing $Y$ we create a curve with a positive $\bar{K}$ label, all our future curves must have a positive $\bar{K}$ label. So since the image of the Hirzebruch curve, once it is of type 1, must have negative $\bar{K}$ label, all $\bar{K}$ labels on $Y$ must be non-positive.

The other restriction to this process comes from the determinant labels. We recall their definition from \cite{DivisorialValuations}. 

\begin{Definition} Suppose $F_i$ and $Y$ are as above. Then the determinant label of the divisorial valuation corresponding to $F_i$ is the determinant of the Gram matrix of minus-intersection form on all curves at infinity of $Y,$ except $F_i$. That is,
$$d_{F_i}= \det \left( -F_j\cdot F_k\right)_{j,k\neq i} $$
\end{Definition}

The determinant labels are much harder to deal with than the $\bar{K}$ labels. In particular, in order to have recursive formulas for them, we need to introduce determinant labels of edges of the graph, and the formulas are somewhat more complicated. Please see \cite{DivisorialValuations} for the details.

If the determinant label of a curve is positive, then one can contract all other curves in the analytic category. If it is negative, then there exists an effective divisor with support in the union of all other curves and positive self-intersection.

\begin{Theorem} Suppose $E$ on $Z$ is a Hirzebruch curve, and it is of type 1. Then $\phi(E)$ must have $\bar{K}$ label (-1) and a positive determinant label.
\end{Theorem}

\proof Suppose $\phi(E)=F.$ By Lemma 2.1, the $\bar{K}$ label of $F$ must divide the $\bar{K}$ label of $E$, so it has to be (-1). 

As for the determinant label of $F$, suppose first that it is negative. Then there is a divisor $D$ supported on all other curves at infinity on $Y$ such that $D^2>0.$ Then $\phi^*(D)^2 >0,$ which contradicts the fact that the determinant label of $E$ is nonnegative.

Suppose now that the determinant label of $F$ is 0. Then $F$ itself must by a Hirzebruch curve on $Y$ (this can be proven similar to Theorem 4.4 of \cite{DivisorialValuations}). One of the fibers of the corresponding fibration consists entirely of curves at infinity. So its full pullback to $Z$ consists of curves at infinity that do not include $E,$ and thus it intersects trivially with $C$, a generic fiber of the Hirzebruch fibration that corresponds to $E$. So $\phi(C)$ intersects trivially with the fiber of the fibration on $Y$, thus it is a fiber itself. This map from a fiber to a fiber is clearly 1-to-1. Therefore by \cite{Gwozdziewicz} $\phi$ is not a Keller map.  
\endproof

\begin{Corollary} All curves on $Y$ have non-negative determinant labels.
\end{Corollary}

\proof One can see (cf. \cite{DivisorialValuations}) that once we create a curve with a negative determinant label, all determinant labels afterwards will be also negative.
\endproof

\begin{Theorem} The curves at infinity on $Y$ generate the Mori cone of effective curves of $Y$.
\end{Theorem}

\proof Suppose $C$ is an irreducible curve on $Y$ that is not a curve at infinity. It is equivalent to $\sum a_iF_i.$ We just need to show that all $a_i$ are nonnegative. Suppose the opposite. Then in the Picard group of $Y$ $C+D_1=D_2,$ where $D_i$ are effective, supported outside of $\mathbb A^2$, have no common support, and $D_1\neq 0$. Multiplying by $D_2,$ we get that $C\cdot D_2 +D_1\cdot D_2=D_2^2.$ So $D_2^2\geq 0,$ with equality if and only if $C\cdot D_2=0$ and $D_1\cdot D_2=0$. Take any curve $F$ in the support of $D_1.$ If $F$ has positive determinant label, we get a contradiction right away. If it has determinant label 0, then one of its branches (the irreducible components of its complement in the graph of curves at infinity) can support a divisor with self-intersection 0. However, $D_1\cdot D_2 =0,$ so $F\cdot D_2=0$. Thus, there is a curve at infinity $F'\neq F$ that intersects positively with $D_2.$ So for a small $\epsilon >0$ $(D_2+\epsilon F')^2 >0,$ a contradiction. 
\endproof

\begin{Corollary} Suppose $C$ is an irreducible curve on $Y$ that is not a curve at infinity. Then either $C^2>0$ or $C^2=0$ and $C$ is a fiber of some Hirzebruch fibration. 
\end{Corollary}

\proof Because $C$ is linearly equivalent to an effective divsor at infinity, $C^2\geq 0.$  If $C^2=0,$ recall that all $\bar{K}$ labels are nonpositive, so $K_Y$ is a linear combination of curves at infinity with negative coefficients. Thus, $C\cdot K_Y<0.$ From the adjunction formula for $C$, this implies that $C$ is rational and smooth, and $C\cdot K_Y=-2$. This means that either $C$ intersects transversally one curve with $\bar{K}$ label (-1), or it intersects transversally two curves with $\bar{K}$ labels 0, or it has intersection of multiplicity 2 with one curve with $\bar{K}$ label 0. In the first case, $C$ is a fiber of a Hirzebruch fibration; the other two cases can be ruled out as follows. From \cite{DivisorialValuations}, the sum of the $\bar{K}$ label and the determinant label of any curve at infinity is always odd. So the determinant labels of 0-curves are strictly positive. The linear combination of curves at infinity that is equivalent to $C$ cannot contain the  curve(s) with $\bar{K}$ label 0 that $C$ intersects, so it must have negative self-intersection, a contradiction.
\endproof

Moreover, similarly to the main idea of \cite{AmpleRamification}, one can prove the following theorem.

\begin{Theorem} Suppose $E$ on $Z$ is a Hirzebruch curve, and it is of type 2 (i.e. $\phi(E)$ is a point). Recall the Stein factorization $\phi=\rho \circ \tau$ with the middle surface $W$. Then for every curve $E_i$ of type 3 on $Z$ either $\tau(E_i)$ contains $\tau(E)$ or $\tau(E_i)$ intersects a curve with $\bar{K}$ label 0 and no other curve at infinity. (Recall that $\tau(E_i)$ has exactly one point at infinity (cf. \cite{AmpleRamification}).
\end{Theorem}

\proof Suppose that $R_i$ are the exceptional curves of type $3$ on $W$, and $r_i\geq 1$ are the corresponding ramification indices. The Keller Map Adjunction Formula asserts that
$$\bar{K}_W=\rho ^* \bar{K}_Y + \sum_{i} r_iR_i$$

Following \cite{AmpleRamification}, we call $R=\sum_{i} r_iR_i$ the di-critical log-ramification divisor.

For each $R_i$ we have the following adjunction inequality (like in \cite{AmpleRamification}, proof of Theorem 3.2):
$$R_i\cdot \bar{K}_W \geq -2 + val(R_i),$$
where $val(R_i),$ the valency of $R_i$, is the number of points on $R_i$, that lie on other curves at infinity on $W.$  As a corollary,
$$\phi_*(R_i)\cdot \bar{K}_Y+R_i\cdot \bar{R}=R_i \cdot (\rho ^* \bar{K}_Y +\bar{R}) =R_i\cdot \bar{K}_W\geq -1$$

If $\rho(R_i)$ intersects on one or more curves with negative $\bar{K}$ labels, $\rho_*(R_i)\cdot \bar{K}_Y \leq -1.$ So
$$R_i\cdot \bar{R} \geq  -1 - \rho_*(R_i) \cdot \bar{K_Y} \geq 0$$

Given that $E$ is of type $2,$ suppose $P=\tau(E).$  Consider  the set $S$ of di-critical curves $R_i$ on $W$ that do not contain $P$ but do intersect some curve of type 1 with negative $\bar{K}$ label. Consider the divisor $D=\sum \limits_{i, R_i\in S}r_iR_i$. Then $D$ does not intersect with any $R_j$ not in its support, so $D^2=D\cdot \bar{R} =\sum (R_i\cdot \bar{R})\geq 0.$ Therefore, $\tau^*D ^2\geq 0$ and $\tau^*D$ does not contain $E$. But this implies that the support of $\tau^*D$ is a union of one or more full fibers of the Hirzebruch fibration corresponding to $E$. If $D\neq 0,$ that is $S$ is not empty, this implies that $\tau^*(D)\cdot E >0,$ so support  of $D$ contains $P$, a contradiction.
\endproof

The proof above also implies that once $\phi(\pi^{-1}_*(\infty))$ is a curve, all images of a type 3 (di-critical) curves can only intersect curves with $\bar{K}$ labels 0. Therefore, during the process of creation of $Y$ the point $\phi(\pi^{-1}_*(\infty))$ initially has all images of the curves of type 3 passing through it, and it is ``losing them" along the way, by creating a curve with $\bar{K}$ label 0 and then going back into the ``negative $\bar{K}$ territory". The following picture shows the curves on $Y$. The red curves are of images of the curves of type 3, the blue dot is  $\phi(\pi^{-1}_*(\infty))$:

\setlength{\unitlength}{0.8 cm}
\begin{picture}(10,5.5)
\thicklines
\put(0,-2){Fig. 32}
\multiput(0,2)(0.5,0){3}{\Huge .}
\put(1.5,1.92){$\circ$}
\put(1.4,1.5){-$1$}
\put(1.5,2.4){$\downarrow$}
\put(0.7,3){$\phi(\pi^{-1}_*(\infty))$ }
\qbezier(0,1)(2,1)(3,0)
\put(3,-0.4){-$1$}
{\color{red}
\qbezier(0,0)(2,0)(3,2)
\qbezier(0,0.5)(2.7,1.65)(3,-1.5)
\put(0.7,-0.5){\small $\phi(type \ 3)$}
}
{\color{blue}
\put(1.67,0.54){$\bullet$}
}

\multiput(5,2)(0.5,0){3}{\Huge .}
\multiput(6.5,1.92)(2,0){2}{$\circ$}
\put(6.71,2.05){\line(1,0){1.84}}
\put(6.4,1.5){-$1$}
\put(8.5,1.5){$0$}
\put(7.5,2.3){$\downarrow$}
\put(6.7,2.9){$\phi(\pi^{-1}_*(\infty))$ }
\qbezier(5,0)(7,0)(8,-1)
\put(8,-1.4){-$1$}
\qbezier(6,-1)(7,0)(9,1)
\put(9,0.5){$0$}
{\color{red}
\qbezier(5,-1)(7,1)(8,-2)
\qbezier(5,1)(8,0)(9,0)
\put(7,-2.3){\small $\phi(type \ 3)$}
\put(8,-0.3){\small $\phi(type \ 3)$}
}

{\color{blue}
\put(6.57,-0.43){$\bullet$}
}

\multiput(11,2)(0.5,0){3}{\Huge .}
\multiput(12.5,1.92)(1,0){3}{$\circ$}
\multiput(12.71,2.05)(1,0){2}{\line(1,0){0.84}}
\put(12.4,1.5){-$1$}
\put(13.4,1.5){-$1$}
\put(14.5,1.5){$0$}
\put(13.1,2.2){$\swarrow$}
\put(13.5,2.3){$\downarrow$}
\put(13.65,2.2){$\searrow$}
\put(12.7,2.8){$\phi(\pi^{-1}_*(\infty))$ }

\put(5.89,3.71){$\searrow$}
\put(5.9,3.71){$\searrow$}
\qbezier(1.7,3.7)(4,6)(6,4)
\put(12.89,3.71){$\searrow$}
\put(12.9,3.71){$\searrow$}
\qbezier(8.7,3.7)(11,6)(13,4)

\end{picture}
\vskip 3cm

It makes sense to stop the process after all curves of type 3 are sent to intersect exclusively with curves with $\bar{K}$ labels 0; this is how our frameworks were obtained.

Another feature of our frameworks has to do with a certain 0-curve on Z. It is the first 0-curve  that was created on Z, the curve number 3 on Figure 7 for the First Framework  and Figure 21 for the Second Framework. It was proven in \cite{DivisorialValuations}, Theorems 4.1 and  4.2, that every curve on $Z$ with a negative $\bar{K}$ and negative determinant label must have such curve as an ancestor. All type 1 curves that are mapped to the original line at infinity on the $\mathbb P^2$ (the curve number 1 on $Y$) have this property. Moreover, they are ``simultaneously determinant-negative": removing them all from the graph of $Z$ produces negative-definite minus-self-intersection form. These collections of curves were discussed in section 5 of \cite{DivisorialValuations}. It can actually be shown that they all must lie  ``outside of a single Hirzebruch curve" (curve number 2  on Figure 7 for the First Framework  and Figure 21 for the Second Framework).

It is important to understand that while some parts of our frameworks may appear random or miraculous, they are not. In particular, the degrees of the maps on the chain of curves between the forked curves with $\bar{K}$ labels (-5) and (-2) (13 and 23, depending in the framework) can be calculated as the index of the cyclic quotient singularity obtained by contracting the chain of the curves on $Y$ between the two forked curves there.  The possible types of maps from the branches also have toric origin. The number of various branches can be calculated from the total ramification of the Belyi maps.

Several questions, some more concrete than the others, naturally appear in relation to our frameworks.

\begin{Question} Is there a simple reason why in the First Framework there is no map $\phi$? If so, it would be really helpful, as it might help pre-screen any further framework examples, before embarking on tedious and time-consuming computer calculations.
\end{Question}

\begin{Question} Can one construct an explicit infinite series of frameworks? This definitely seems possible, even with just two Belyi maps.
\end{Question}

\begin{Question} Can one formalize the notion of a framework, and to actually find ALL solutions to this combinatorial problem (or, more realistically, all ``small" solutions, in some reasonable sense)?
\end{Question}

\begin{Question} There are several questions regarding the notion of the isotope. Are there any more isotopes of the First Framework? Do the Second Framework and the Three-dessin Framework have other isotopes, and can they be classified? Does every framework have only finitely many isotopes?
\end{Question}

\begin{Question} Do our frameworks actually provide maps from the tubular neighborhood of the union of curves at infinity on $Z$, without the curves of type 3 and 4, to the tubular neighborhood of the union of curves at infinity on $Y$? It seems like we get maps between the abelianizations of the fundamental groups, but is it enough?
\end{Question}

\begin{Question} The smallest topological degree of the Keller maps that would come from our frameworks is 16. The current best lower bound for a topological degree of a Keller map is 6 (\cite{Zoladek}). Can one use some ideas from Section 6 to greatly improve this bound?
\end{Question}

\begin{Question} (The biggest question of all). Can one actually use our frameworks to contruct a Keller map? If you have some time and knowledge in computing, I am very open to collaboration, and will be glad to share with you many further details beyond the discussion at the end of section 3.
\end{Question}

{\bf Acknowledgments.} The author is forever indebted to Vasilii Alekseevich Iskovskikh, who introduced him, and many others, to the beauty of birational geometry. The author is also indebted to David Wright, Ed Formanek, and Adrian Vasiu for stimulating discussions related to the Jacobian Conjecture. Additionally, the author wants to thank his colleagues at Binghamton University and University of Pittsburgh, in particular Marcin Mazur, Bogdan Ion, and Thomas Hales, for their interest and patient support.

\end{document}